\documentclass[aos,dvips]{arximspdf}
\usepackage{mathbh}
\usepackage{mathrsfs}
\usepackage{graphicx}

\doi{10.1214/07-AOS527}
\volume{36}
\issue{4}
\pubyear{2008}
\firstpage{1925}
\lastpage{1956}

\makeatletter
\newtheorem{theorem}{Theorem}
\newproclaim{nation}{Notational convention}
\newtheorem{lemma}[theorem]{Lemma}
\newproclaim{definition}{Definition}
\newtheorem{proposition}[theorem]{Proposition}
\newtheorem{corollary}[theorem]{Corollary}
\newproclaim{remark}{Remark}
\newproclaim{condition}{Condition}
\def\allWA{\textup{(\hyperref[item:Wreg]{W-1})}--\textup{(\hyperref[item:Wvstd]{W-5})}}
\newcommand{\eqref}[1]{(\ref{#1})}
\def\K{\sigma^2}
\newcommand\nj[1]{{n}_{#1}}
\def\maxscale{J}
\def\nmaxscale{\maxscale_{n}}
\def\lowscale{L}
\def\nlowscale{\lowscale_{n}}
\def\upscale{U}
\def\nupscale{\upscale_{n}}
\def\rme{e}
\def\empSbis{\tilde{\mathbf{S}}}
\def\Sclt{\widehat{{\mathbf{S}}}}
\def\Scltbar{\overline{\mathbf{S}}}
\newcommand\jmean[1]{\langle #1 \rangle}
\newcommand\cardinal[1]{|#1|}
\def\admin{d_{\min}}
\def\admax{d_{\max}}
\def\dj{u}
\def\rmi{i}
\def\cl{\mathop{\stackrel{\mathcal{L}}{\longrightarrow}}}
\def\vareta{\mathcal{J}}
\def\hd{\hat{d}} 
\def\hs{\hat{\sigma}} 
\def\td{d_{0}} 
\def\htd{\tilde{d}} 
\def\dwt{W}
\def\bdwt{\mathbf{\dwt}}
\def\indexset{\mathcal{I}}
\def\LWF{LWF}
\def\ELWF{ELWF}
\def\L{\mathrm{T}}
\def\be{\mathbf{e}}
\def\bmu{\boldsymbol{\mu}}
\def\diffop{\mathbf{\Delta}}
\def\1{\mathbh{1}}
\def\Nset{\mathbb{N}} 
\def\Zset{\mathbb{Z}}
\def\Rset{\mathbb{R}} 
\def\PE{\mathbb{E}} 
\def\PVar{\operatorname{Var}}
\def\PCov{\operatorname{Cov}}
\def\Cum{\operatorname{Cum}}
\def\calN{\mathcal{N}}
\def\calB{\mathcal{B}}
\def\calH{\mathcal{H}}
\def\ie{that is, }
\def\eg{e.g., }
\def\argmin{\mathrm{Argmin}}
\def\prob{\mathbb{P}}
\newcommand{\eqdef}{\mathop{\stackrel{\mathrm{def}}{=}}}
\def\contrastsymb{\mathrm{L}}
\def\contrastL{\hat{\contrastsymb}}
\def\contrast{\tilde{\contrastsymb}}
\def\contrastdet{\contrastsymb}
\def\fluctuation{\mathrm{E}}
\newcommand{\AsympVarWWE}[2][]
{\mathrm{V}(#2)}
\def\vjsymb{\sigma}
\newcommand{\hvj}[3][]{
\ifthenelse{\equal{#1}{}}{\hat{\vjsymb}^2_{#2}}{\hat{
\vjsymb}^2_{#2}(#1)}}
\newcommand{\Kvar}[2][]{
\mathrm{K}(#2)}
\newcommand{\Kvarsq}[2][]{
\mathrm{K}^2(#2)}
\def\mes{\nu}
\def\d{d}
\def\densraid{\mathbf{D}}
\newcommand{\bdens}[4][]{\densraid_{#2}({#3};#4)}
\def\densasympletter{D}
\newcommand{\bdensasymp}[4][]{\mathbf{\densasympletter}_{\infty,#2}({#3};#4)}
\def\intbdensletter{I}
\newcommand{\intbdens}[3][]{\mathrm{\intbdensletter}_{#2}(#3)}
\renewcommand{\k}{\ell}
\makeatother

\begin{document}
\begin{frontmatter}

\title{A wavelet whittle estimator of the memory parameter of a
nonstationary\break Gaussian time series\protect\thanksref[1]{T1}}
\runtitle{A wavelet whittle estimator}

\begin{aug}
\author[A]{\fnms{E.} \snm{Moulines}\ead[label=e1]{moulines@tsi.enst.fr}},
\author[A]{\fnms{F.} \snm{Roueff}\corref{}\ead[label=e2]{roueff@tsi.enst.fr}}
\and
\author[B]{\fnms{M. S.} \snm{Taqqu}\ead[label=e3]{murad@math.bu.edu}}
\runauthor{E. Moulines, F. Roueff and M. S. Taqqu}
\thankstext[1]{T1}{Supported in part by NSF Grant DMS-05-05747.}
\affiliation{T\'{e}l\'{e}com Paris/CNRS LTCI and Boston University}
\address[A]{E. Moulines\\
F. Roueff\\
E.N.S.T. 46\\
rue Barrault\\
75634 Paris C\'{e}dex 13\\
France\\
\printead{e1}\\
\phantom{\textsc{E-mail:\ }}\printead*{e2}} 
\address[B]{M. S. Taqqu\\
Department of Mathematics and Statistics\\
Boston University\\
Boston, Massachusetts 02215\\
USA\\
\printead{e3}}
\end{aug}

\received{\smonth{7} \syear{2007}}
\revised{\smonth{7} \syear{2007}}

%
\begin{abstract}
We consider a time series $X=\{X_k,\,k\in{\mathbb{Z}}\}$ with memory parameter
$d_{0}\in\mathbb{R}$. This time series is
either stationary or can be made stationary after differencing a finite
number of times. We study the ``local Whittle
wavelet estimator'' of the memory parameter $d_{0}$. This is a
wavelet-based semiparametric
pseudo-likelihood maximum method estimator.
The estimator may depend on a given finite range of scales
or on a range which becomes infinite with the sample size. We show
that the estimator is consistent and rate
optimal if $X$ is a linear process, and is asymptotically normal if $X$
is Gaussian.
\end{abstract}

%

\end{frontmatter}
%

\section{Introduction}

Let $X\eqdef\{X_k\}_{k\in\Zset}$ be a process, not necessarily
stationary or invertible. Denote by $\diffop
X$, the first order difference, $(\diffop X)_\k=X_{\k}-X_{\k-1}$,
and by $\diffop^k X$, the $k$th order difference.
Following \cite{hurvichray1995}, the process $X$ is said to have
memory parameter $\td$, $\td\in\Rset$,
if for any integer $k>\td-1/2$, $U\eqdef\diffop^{k}X$ is covariance
stationary with spectral measure
%
\begin{equation}\label{eq:spdelta}
\mes_{U}(d\lambda)=|1-\rme^{-\rmi\lambda}|^{2(k-\td)}\mes^\ast
(d\lambda),\qquad\lambda\in[-\pi,\pi],
\end{equation}
where $\mes^\ast$ is a nonnegative symmetric measure on $[-\pi,\pi]$
such that, in a neighborhood of the origin,
it admits a positive and bounded density. The process $X$ is covariance
stationary if and only if $\td<1/2$. When $\td>0$,
$X$ is said to exhibit long memory or long-range dependence. The
\textit{generalized spectral measure} of $X$ is
defined as
%
\begin{equation}\label{eq:fmodele}
\mes(\d\lambda)\eqdef|1-\rme^{-\rmi\lambda}|^{-2\td}\mes
^\ast(\d
\lambda),\qquad\lambda\in[-\pi,\pi].
\end{equation}

We suppose that we observe $X_1,\dots,X_n$ and want to
estimate the exponent $\td$ under the following semiparametric set-up
introduced in \cite{robinson1995g}.
Let $\beta\in(0,2]$, $\gamma>0$ and $\varepsilon\in(0,\pi]$, and
assume that
\[
\mes^\ast\in\calH(\beta,\gamma,\varepsilon),
\]
where $\calH(\beta,\gamma,\varepsilon)$ is the class of finite
nonnegative symmetric measures on $[-\pi,\pi]$
whose restrictions on $[-\varepsilon,\varepsilon]$ admit a density
$g$, such that, for all
$\lambda\in(-\varepsilon,\varepsilon)$,
%
\begin{equation}\label{eq:Hbeta}
|g(\lambda)- g(0) | \leq\gamma g(0)  |\lambda|^\beta.
\end{equation}

Since $\varepsilon\leq\pi$, $\mes^\ast\in\calH(\beta,\gamma
,\varepsilon
)$ is only a local condition for $\lambda$ near 0. For instance, $\mes
^\ast$
may contain atoms at frequencies in $(\varepsilon,\pi]$ or have an
unbounded density on this domain.

We shall estimate $\td$ using the semiparametric \textit{local Whittle
wavelet estimator} defined in
Section~\ref{sec:wavel-whittle-estim}. We will show that under
suitable conditions, this estimator is consistent
(Theorem~\ref{thm:Rates}), the convergence rate is optimal
(Corollary~\ref{coro:optimal-rate}) and it is asymptotically
normal (Theorem~\ref{thm:CLT}). In Section~\ref{sec:disc-concl-remarks}, we discuss how it compares to other
estimators.

There are two popular semiparametric estimators for the memory
parameter $\td$ in the frequency domain:

\begin{enumerate}[(1)]
\item[(1)] the Geweke--Porter--Hudak (GPH) estimator introduced in \cite
{gewekeporterhudak1983} and analyzed
in \cite{robinson1995l}, which involves a
regression of the log-pe\-rio\-do\-gram on the log of low frequencies;

\item[(2)] the local Whittle (Fourier) estimator (or \LWF) proposed
in \cite{kunsch1987} and developed in \cite{robinson1995g}, which is
based on the Whittle approximation of the Gaussian likelihood,
restricted to low
frequencies.
\end{enumerate}
Corresponding approaches may be considered in the wavelet domain. By
far, the most widely used wavelet estimator is based on
the log-regression of the wavelet coefficient variance on the scale
index, which was introduced in
\cite{abryveitch1998}; see also~\cite{moulinesrouefftaqqu2007a}
and~\cite{moulinesrouefftaqqu2006b} for recent
developments.
A wavelet analog of the \LWF, referred to as the \textit{local Whittle
wavelet estimator} can also be defined.
This estimator was proposed for analyzing noisy data in a parametric context in \cite{wornelloppenheim1992}
and was considered by several authors, essentially in a parametric
context (see, \eg\cite{kaplankuo1993} and
\cite{mccoywalden1996}).
To our knowledge, its theoretical properties are not known (see the
concluding remarks in \cite{velasco1999}, page~107).
The main goal of this paper is to fill this gap in a semiparametric context.
The paper is structured as follows. In Section~\ref
{sec:discr-wavel-analys}, the
wavelet analysis of a time series is presented and some results on the
dependence structure of the wavelet coefficients are
given. The definition and the asymptotic properties of the local
Whittle wavelet estimator are given in
Section~\ref{sec:wavel-whittle-estim}: the
estimator is shown to be rate optimal under a general condition on the
wavelet coefficients, which are satisfied when $X$ is
a linear process with four finite moments, and it is shown to be
asymptotically normal under the additional condition that $X$
is Gaussian. These results are discussed in Section~\ref
{sec:disc-concl-remarks}. The proofs can be found in the remaining
sections. The linear case is considered in Section~\ref
{sec:assumpcons}. The asymptotic
behavior of the wavelet Whittle likelihood is studied in
Section~\ref{sec:asympt-behav-contr} and weak consistency is studied in
Section~\ref{sec:weak-consistency}. The proofs of the main results
are gathered in Section~\ref{sec:rates}.

\section{The wavelet analysis}\label{sec:discr-wavel-analys}

The functions $\phi(t)$, $t\in\Rset$, and $\psi(t)$, $t\in\Rset$, will
denote the father and mother wavelets respectively,
and
$\hat{\phi}(\xi)\eqdef\int_{\Rset} \phi(t)\rme^{-\rmi\xi t}\,
\d t$ and
$\hat{\psi}(\xi)\eqdef\int_{\Rset}
\psi(t)\rme^{-\rmi\xi t}\,\d t$ their Fourier
transforms. We suppose that $\phi$ and $\psi$ satisfy the following
assumptions:

\begin{longlist}[(W-5)]

\item[(W-1)]\label{item:Wreg} $\phi$ and $\psi$ are integrable and have
compact supports, $\hat{\phi}(0) =\break \int_{\Rset} \phi(x)\, \d x =
1$ and $\int_{\Rset} \psi^2(x)\, \d x = 1$;

\item[(W-2)]\label{item:psiHat}
there exists $\alpha>1$ such that
$\sup_{\xi\in\Rset}|\hat{\psi}(\xi)|(1+|\xi|)^{\alpha}
<\infty$;

\item[(W-3)]\label{item:MVM} the function $\psi$ has $M$ vanishing moments,
\ie $ \int_{\Rset} t^l \psi(t) \,\d t=0$ for all $l=0,\dots,M-1$;

\item[(W-4)]\label{item:MIM} the function $ \sum_{k\in\Zset} k^l\phi
(\cdot-k)$ is a polynomial of degree $l$ for all $l=0,\dots,M-1$;

\item[(W-5)]\label{item:Wvstd} $\td$, $M$, $\alpha$ and $\beta$ are such that
$(1+\beta)/2-\alpha<\td\leq M$.
\end{longlist}

Assumption~(\hyperref[item:Wreg]{W-1}) implies that $\hat{\phi}$ and $\hat{\psi}$
are everywhere infinitely differentiable.
Assumption~(\hyperref[item:psiHat]{W-2}) is regarded as a \textit{regularity
condition} and assumptions~(\hyperref[item:MVM]{W-3}) and~(\hyperref[item:MIM]{W-4})
are often referred to as \textit{admissibility conditions}.
When~(\hyperref[item:Wreg]{W-1}) holds, assumptions~(\hyperref[item:MVM]{W-3}) and~(\hyperref
[item:MIM]{W-4}) can be expressed in different ways. (\hyperref[item:MVM]{W-3})
is equivalent to asserting that the first $M-1$ derivative of $\hat
{\psi}$ vanish at
the origin and hence
%
\begin{equation}
\label{eq:MVM}
|\hat{\psi}(\lambda)|=O(|\lambda|^{M})\qquad\mbox{as }\lambda
\to0.
\end{equation}
And, by \cite{cohen2003}, Theorem~2.8.1, page~90, (\hyperref[item:MIM]{W-4}) is
equivalent to
%
\begin{equation}
\label{eq:MIM}
\sup_{k\neq0} |\hat{\phi}(\lambda+2k\pi)|=O(|\lambda|^{M})\qquad
\mbox{as }\lambda\to0.
\end{equation}
Finally,~(\hyperref[item:Wvstd]{W-5}) is the constraint on $M$
and $\alpha$ that we will impose on the wavelet-based estimator of the
memory parameter $\td$
of a process having generalized spectral
measure~(\ref{eq:fmodele}) with $\mes^\ast\in\calH(\beta,\gamma
,\varepsilon)$ for some positive $\beta$, $\gamma$ and $\varepsilon$.
Remarks~\ref{rem:DefBdensAsymp} and~\ref{rem:W5} below provide some
insights into (\hyperref[item:Wvstd]{W-5}).
We may consider nonstationary processes $X$ because the wavelet analysis
performs an implicit differentiation of order $M$. It is perhaps less
well known that, in addition, wavelets can be used
with noninvertible processes ($\td\leq-1/2$) due to the regularity
condition~(\hyperref[item:psiHat]{W-2}). These two properties of
the wavelet are, to some extent, similar to the properties of the
tapers used in Fourier analysis
(see, e.g., \cite{hurvichray1995,velasco1999}).

Adopting the engineering convention that large values of the scale
index $j$ correspond to coarse scales (low frequencies),
we define the family $\{\psi_{j,k}, j \in\Zset, k \in\Zset\}$ of
translated and dilated functions,
$\psi_{j,k}(t)=2^{-j/2}\psi(2^{-j}t-k)$, $j\in\Zset$, $k\in\Zset$.
If $\phi$ and $\psi$ are the scaling and wavelet functions
associated with a multiresolution analysis
(see \cite{cohen2003}), then $\{\psi_{j,k},j\in\Zset,k\in\Zset
\}$
forms an orthogonal basis in $L^2(\Rset)$.
A standard choice are
the Daubechies wavelets (DB-$M$), which are parameterized by the number
of their vanishing moments
$M$. The associated scaling and wavelet functions $\phi$ and $\psi$
satisfy (\hyperref[item:Wreg]{W-1})--(\hyperref[item:MIM]{W-4}),
where $\alpha$ in~(\hyperref[item:psiHat]{W-2}) is a function of $M$ which
increases to infinity as $M$ tends to infinity
(see \cite{cohen2003}, Theorem~2.10.1). In this work, however, we
neither assume that
the pair $\{\phi,\psi\}$ is associated with a multiresolution analysis
(MRA), nor that the $\psi_{j,k}$'s form
a Riesz basis. Other possible choices are discussed
in~\cite{moulinesrouefftaqqu2007a},
Section~3.

The wavelet coefficients of the process $X =\{X_\k, \k\in\Zset\}$ are
defined by
%
\begin{equation}\label{eq:coeff}
\dwt_{j,k}\eqdef\int_{\Rset} X(t) \psi_{j,k}(t)\,\d t,\qquad j \geq
0, k
\in\Zset,
\end{equation}
where $X(t) \eqdef\sum_{k\in\Zset} X_k \phi(t-k)$. If $(\phi
,\psi)$
define an MRA, then $X_k$ is identified with the
$k$th approximation coefficient at scale $j=0$ and $\dwt_{j,k}$ are
the details coefficients at scale $j$.

Because translating the functions $\phi$ or $\psi$ by an integer
amounts to translating the sequence $\{\dwt_{j,k},k\in\Zset\}$ by
the same integer for all $j$, we can suppose, without loss of
generality,
that the supports of $\phi$ and $\psi$ are included in
$[-\L,0]$ and $[0,\L]$, respectively, for some integer $\L\geq1$.
Using this convention, it is easily seen that the wavelet coefficient
$\dwt_{j,k}$ depends only on the available observations
$\{X_1, \dots, X_n\}$ when $j\geq0$ and $0\leq k < \nj{j}$, where,
denoting the integer part of $x$ by $[x]$,
%
\begin{equation}\label{eq:Nj}
\nj{j}\eqdef\max\bigl([2^{-j}(n-\L+1)-\L+1], 0\bigr) .
\end{equation}

Suppose that $X$ is a (possibly nonstationary) process with memory
parameter $\td$ and generalized spectral measure
$\mes$. If $M>\td-1/2$, then $\diffop^M X$ is stationary and hence, by
\cite{moulinesrouefftaqqu2007a}, Proposition~1,
the sequence of wavelet coefficients $\dwt_{j,\centerdot}$ is a
stationary process and we can define $\sigma^{2}_{j}(\nu) \eqdef
\PVar(\dwt_{j,k})$.
Our estimator takes advantage of the \textit{scaling} and \textit{weak
dependence} properties of the wavelet coefficients,
as expressed in the following condition, which will be shown to hold in
many cases of interest.
\begin{condition}
\label{cond:GeneralCondition}
There exist $\beta> 0$ and $\K>0$ such that
%
\begin{equation}\label{eq:BiasTermAssump}
\sup_{j\geq1} 2^{\beta j} \biggl|\frac{\sigma^{2}_{j}(\nu)}{\K 2^{2
\td j}}-1\biggr| < \infty
\end{equation}
and
%
\begin{equation}
\label{eq:rosenthalAssump}
\sup_{n \geq1} \sup_{j=1,\dots,\nmaxscale}  (1+ \nj{j}
2^{-2j\beta
})^{-1}
\nj{j}^{-1} \PVar\Biggl(\sum_{k=0}^{\nj{j}-1} \frac{\dwt
_{j,k}^{2}}{\sigma^{2}_{j}(\nu)}\Biggr) < \infty .
\end{equation}
\end{condition}

Equation (\ref{eq:BiasTermAssump}) states that, up to the multiplicative
constant $\K$,
the variance $\sigma^{2}_{j}(\nu)$ is approximated by $2^{2\td j}$ and that
the error goes to zero exponentially fast as
a function of $j$. It is a direct consequence of the approximation of
the covariance of the wavelet coefficients established in~\cite
{moulinesrouefftaqqu2007a}.
Equation (\ref{eq:rosenthalAssump}) imposes a bound on the variance of the
normalized partial sum of the stationary
centered sequence $\{ \sigma^{-2}_{j}(\nu) \dwt_{j,k}^2 \}$, which, provided
that $\nj{j}2^{-2j\beta}=O(1)$, is equivalent to what occurs when
these variables are independent. We stress that the wavelet
coefficients $\dwt_{j,k}$ are, however,
\textbf{not} independent, nor can they be approximated by independent
coefficients; see~\cite{moulinesrouefftaqqu2007a}.
Establishing~(\ref{eq:rosenthalAssump}) requires additional assumptions
on the process $X$ that go beyond its covariance
structure since $\dwt_{j,k}^2$ is involved; see Theorem~\ref{thm:AssumpCons}, where this
property is established for a general class of linear processes.
We have isolated relations~(\ref{eq:BiasTermAssump}) and~(\ref
{eq:rosenthalAssump}) because in our semiparametric context,
these two relations are sufficient to show that
the wavelet Whittle estimator converges to $\td$ at the optimal rate
(see Theorem~\ref{thm:Rates} below).\looseness=1

Let us recall some definitions and results from~\cite
{moulinesrouefftaqqu2007a} which are used here.
As noted above, for a given scale $j$, the process $\{ \dwt_{j,k} \}_{k
\in\Zset}$ is covariance stationary. It will be
called the \textit{within-scale} process because all the $\dwt_{j,k}$, $k
\in\Zset$, share the same~$j$.
The situation is more complicated when considering two different scales
$j > j'$ because
the two-dimensional sequence $\{[\dwt_{j,k},\dwt_{j',k}]^T \}_{k\in
\Zset}$ is not
stationary, as a consequence of the pyramidal wavelet scheme.
A convenient way to define a joint spectral density for wavelet
coefficients is to consider the \textit{between-scale}
process.
\begin{definition}
The sequence $\{[\dwt_{j,k},\bdwt_{j,k}(j-j')^T]^T \}_{k\in\Zset
}$, where
\[
\bdwt_{j,k}(j-j') \eqdef[\dwt_{j',2^{j-j'}k},\dots, \dwt
_{j',2^{j-j'}k+2^{j-j'}-1}]^T,
\]
is called the \textit{between-scale} process at scales $0\leq j'\leq j$.
$\bdwt_{j,k}(j-j')$ is a $2^{j-j'}$-dimensional vector
of wavelet coefficients at scale $j'$.
\end{definition}

Assuming that the generalized spectral measure of $X$ is given by
\eqref{eq:fmodele} and
provided that $M>\td-1/2$, since $\diffop^MX$ is stationary, both the
within-scale process and the
between-scale process are covariance stationary; see~\cite
{moulinesrouefftaqqu2007a}.
Let us consider the case $\mes^\ast\in\calH(\beta,\gamma,\pi)$, that
is, $\varepsilon=\pi$, so that $\mes^\ast$ admits a density
$f^\ast$ in the space $\calH(\beta,\gamma)$ as defined in~\cite
{moulinesrouefftaqqu2007a} and $\mes$ admits a density
$f(\lambda)\eqdef|1-\rme^{-\rmi\lambda}|^{-2\td}f^\ast(\lambda)$.
We denote by $\bdens[\phi,\psi]{j,0}{\cdot}{f}$ the spectral
density of
the within-scale process at scale index $j$ and
by $\bdens[\phi,\psi]{j,j-j'}{\cdot}{f}$ the cross spectral density
between $\{\dwt_{j,k}\}_{k\in\Zset}$
and $\{\bdwt_{j,k}(j-j') \}_{k\in\Zset}$ for $j'<j$. It will be
convenient to set $\dj=j-j'$.
Theorem~1 in~\cite{moulinesrouefftaqqu2007a} states that, under
\allWA,
for all $\dj\geq0$, there exists $C>0$ such that for all $\lambda
\in
(-\pi,\pi)$ and
$j\geq\dj\geq0$,
%
\begin{equation}\label{eq:DjApprox}
|\bdens[\phi,\psi]{j,\dj}{\lambda}{f} - f^\ast(0)
\bdensasymp
[\psi]{\dj}{\lambda}{\td}  2^{2j \td}|
\leq C f^\ast(0)  2^{(2\td-\beta)j},
\end{equation}
where, for all $\dj\geq0$, $d\in(1/2-\alpha,M]$ and $\lambda\in
(-\pi,\pi)$,
%
\begin{equation}\label{eq:bDpsi}
\hspace*{26pt}\bdensasymp[\psi]{\dj}{\lambda}{d} \eqdef
\sum_{l\in\Zset} |\lambda+2l\pi|^{-2d}\be_{\dj}(\lambda+2l\pi)
\overline{\hat{\psi}(\lambda+2l\pi)}\hat{\psi}\bigl(2^{-\dj}(\lambda
+2l\pi)\bigr),
\end{equation}
with $\be_\dj(\xi) \eqdef2^{-u/2} [1, \rme^{-\rmi2^{-\dj}\xi},
\dots,
\rme^{-\rmi(2^{\dj}-1)2^{-\dj}\xi}]^T $.
\begin{remark}\label{rem:DefBdensAsymp}
The condition~(\hyperref[item:Wvstd]{W-5}) involves an upper and a lower bound. The
lower bound guarantees that the series
defined by the right-hand side of~(\ref{eq:bDpsi}) omitting the term
$l=0$ converges uniformly for $\lambda\in(\pi,\pi)$. The
upper bound guarantees that the term $l=0$ is bounded at $\lambda=0$.
As a result, $\bdensasymp[\psi]{\dj}{\lambda}{d}$ is
bounded on $\lambda\in(\pi,\pi)$ and, by~(\ref{eq:DjApprox}), so is
$\bdens[\phi,\psi]{j,\dj}{\lambda}{f}$. In particular, the
wavelet coefficients are short-range dependent. For details, see
the proof of Theorem~1
in~\cite{moulinesrouefftaqqu2007a}.
\end{remark}
\begin{remark}\label{rem:vaepsileqpi}
We stress that~(\ref{eq:DjApprox}) may no longer hold if we only assume
$\mes^\ast\in\calH(\beta,\gamma,\varepsilon)$ with
$\varepsilon<\pi$ since in this case, no condition is imposed on
$\mes
(d\lambda)$ for $|\lambda|>\varepsilon$ and hence
$\dwt_{j,\centerdot}$ may not have a density for all $j$. However, this
difficulty can be circumvented by decomposing
$\mes^\ast$ as
%
\begin{equation}
\label{eq:decompmesstar}
\mes^\ast(\d\lambda)=f^\ast(\lambda)\,\d\lambda+\tilde{\mes
}^\ast(\d
\lambda),
\end{equation}
where $f^\ast$ has support in $[-\varepsilon,\varepsilon]$ and
$\tilde
{\mes}^\ast([-\varepsilon,\varepsilon])=0$; see the
proof of Theorem~\ref{thm:AssumpCons}.
\end{remark}

Here is a simple interpretation of the bound~(\ref{eq:DjApprox}).
For any $d\in\Rset$, $ 2^{2j d}\bdensasymp[\psi]{\dj}{\cdot}{d}$ is
the spectral density of the
wavelet coefficient of the generalized fractional Brownian motion (GFBM)
$\{B_{(d)}(\theta)\}$ defined as the Gaussian process indexed by test
functions $ \theta\in\Theta_{(d)} = \{
\theta\dvtx \int_\Rset|\xi|^{-2d}|\hat{\theta}(\xi)|^2\,\d\xi<
\infty
\}$
with mean zero and covariance
%
\begin{equation}\label{eq:CovFBMGen}
\PCov\bigl(B_{(d)}(\theta_1),B_{(d)}(\theta_2)\bigr) = \int
_{\Rset}
|\xi|^{-2d}\hat{\theta_1}(\xi)\overline{\hat{\theta_2}(\xi
)}\,\d\xi.
\end{equation}
When $d>1/2$, the condition $\int|\xi|^{-2d}|\hat{\theta}(\xi
)|^2\,d\xi< \infty$ requires that $\hat{\theta}(\xi)$ decays sufficiently
quickly at the origin and when $d<0$, it requires that $\hat{\theta
}(\xi)$ decreases sufficiently rapidly at
infinity. Provided that $d\in(1/2-\alpha,M+1/2)$, the wavelet function
$\psi$ and its scaled and translated versions $\psi_{j,k}$ all belong
to $\Theta_{(d)}$. Defining the discrete wavelet transform of
$B_{(d)}$ as $\dwt_{j,k}^{(d)} \eqdef B_{(d)}(\psi_{j,k}), j\in
\Zset
,k\in\Zset$ and $\bdwt^{(d)}_{j,k}(\dj)\eqdef[\dwt
^{(d)}_{j-\dj
,2^{\dj}k},\dots, \dwt^{(d)}_{j-\dj,2^{\dj}k+2^{\dj}-1}]$,
one obtains
%
\begin{equation}\label{eq:CovFBMGenPsijk}
\PCov\bigl(\dwt_{j,k}^{(d)},\bdwt^{(d)}_{j,k'}(\dj)
\bigr)=2^{2dj}
\int_{-\pi}^\pi\bdensasymp[\psi]{\dj}{\lambda}{d}
\rme^{\rmi\lambda(k-k')}\,\d\lambda ;
\end{equation}
see \cite{moulinesrouefftaqqu2007a}, Remark~5, for more details.
Equation (\ref{eq:DjApprox}) shows that the within- and between-scale
spectral densities
$\bdens[\phi,\psi]{j,\dj}{\lambda}{\mes}$ of the process $X$ with
memory parameter $d$ may be approximated by the corresponding
densities of the wavelet coefficients of the GFBM $B_{(d)}$, with an
$L^\infty$-error bounded by $O(2^{(2\td-\beta)j})$.

The approximation~(\ref{eq:DjApprox}) is a crucial step for proving
that Condition~\ref{cond:GeneralCondition} holds
for linear processes. The following theorem is proved in Section~\ref{sec:assumpcons}.
\begin{theorem}\label{thm:AssumpCons}
Let $X$ be a process having generalized spectral measure \textup{(\ref{eq:fmodele})}
with $\td\in\Rset$ and with
$\mes^\ast\in\calH(\beta,\gamma,\varepsilon)$ such that $f^\ast
(0)\eqdef\d\mes^\ast/\d\lambda_{|\lambda=0}>0$, where
$\gamma>0$, $\beta\in(0,2]$ and $\varepsilon\in(0,\pi]$. Then, under
\allWA, the bound \textup{(\ref{eq:BiasTermAssump})} holds with
$\K=f^\ast(0) \Kvar[\psi]{\td}$, where
%
\begin{equation}\label{eq:Kpsi}
\Kvar[\psi]{d}\eqdef\int_{-\infty}^{\infty} |\xi|^{-2d}|\hat
\psi(\xi
)|^2\,\d\xi\qquad\mbox{for any $d\in(1/2-\alpha,M+1/2)$}.
\end{equation}
Suppose, in addition, that there exist an integer $k_0\leq M$ and a
real-valued sequence $\{a_k\}_{k\in\Zset} \in\ell^2(\Zset)$ such that
%
\begin{equation}
\label{eq:linearProc}
(\diffop^{k_0}X)_k=\sum_{t\in\Zset} a_{k-t} Z_t,\qquad k\in\Zset,
\end{equation}
where $\{Z_t\}_{t \in\Zset}$ is a weak white noise process such that
$\PE[Z_t]=0$, $\PE[Z_t^2]=1$, $\PE[Z_t^4]=\PE[Z_1^4]<\infty$ for
all $t\in\Zset$ and
%
\begin{equation}
\label{eq:cumFour}
\Cum(Z_{t_1},Z_{t_2},Z_{t_3},Z_{t_4})=\cases{
\PE[Z_1^4]-3, &\quad $\mbox{if }t_1=t_2=t_3=t_4$,\cr
0, &\quad otherwise.}
\end{equation}
Then, under \allWA, the bound \textup{(\ref{eq:rosenthalAssump})} holds and
Condition~\ref{cond:GeneralCondition} is satisfied.
\end{theorem}
\begin{remark}
Relation~(\ref{eq:rosenthalAssump}) does not hold for every long-memory
process $X$, even with
arbitrary moment conditions; see~\cite{fayroueffsoulier2007}.
\end{remark}
\begin{remark}
Any martingale increment process with constant finite fourth moment, as
in the assumption A3$'$ considered
in~\cite{robinson1995g}, satisfies~(\ref{eq:cumFour}). Another
particular case is given by the following corollary, proved
in Section~\ref{sec:assumpcons}.
\end{remark}

The following result specializes Theorem~\ref{thm:AssumpCons} to a
Gaussian process $X$ and shows that at large scales, the
wavelet coefficients of $X$ can be approximated by those of a process
$\bar{X}$ whose spectral measure $\bar{\mes}$ satisfies
the global condition $\bar{\mes}\in\calH(\beta,\gamma,\pi)$.
\begin{corollary}\label{cor:C1gaussian}
Let $X$ be a Gaussian process having generalized spectral measure~\textup{(\ref{eq:fmodele})}
with $\td\in\Rset$ and with $\mes^\ast\in\calH(\beta,\gamma
,\varepsilon)$ such that
$f^\ast(0)\eqdef\d\mes^\ast/\break\d\lambda_{|\lambda=0}>0$, where
$\gamma
>0$, $\beta\in(0,2]$ and $\varepsilon\in(0,\pi]$.
Then, under \allWA, Condition~\ref{cond:GeneralCondition} is
satisfied with
$\K=f^\ast(0) \Kvar[\psi]{\td}$.

There exists, moreover, a Gaussian process $\overline{X}$ defined on
the same probability space as $X$ with generalized spectral measure
$\bar{\mes}\in\calH(\beta,\gamma,\pi)$ and wavelet coefficients
$\{
\overline{W}_{j,k}\}$ such that
%
\begin{eqnarray}\label{eq:gaussianCaseApproxXbar}
&&\sup_{n\geq1,j\geq0}
\bigl\{\nj{j}2^{j(1+2\td-2\alpha)}+\nj{j}^22^{2j(1-2\alpha
)}\bigr\}^{-1}\nonumber\\[-8pt]\\[-8pt]
&&\qquad{}\times\PE\Biggl[\Biggl|\sum_{k=0}^{\nj{j}-1} \dwt_{j,k}^{2}- \sum
_{k=0}^{\nj
{j}-1} \overline{\dwt}_{j,k}^{2}\Biggr|^2\Biggr]
< \infty .\nonumber
\end{eqnarray}
\end{corollary}

\section{Asymptotic behavior of the local Whittle wavelet
estimator}\label{sec:wavel-whittle-estim}
We first define the estimator. Let $\{c_{j,k},(j,k)\in\indexset\}$ be
an array of centered independent Gaussian random variables
with variance $\PVar(c_{j,k})=\vjsymb_{j,k}^2$, where $\indexset$ is a
finite set. The
negative of its log-likelihood is
$ (1/2) \sum_{(j,k) \in\indexset} \{ c_{j,k}^2/\vjsymb
_{j,k}^2 +
\log(\vjsymb_{j,k}^2)\}$, up to a constant additive term.
Our local Whittle wavelet estimator (LWWE) uses such a contrast process
to estimate the memory parameter $\td$ by
choosing $c_{j,k}=\dwt_{j,k}$. The scaling and weak dependence in
Condition~\ref{cond:GeneralCondition} then suggest
the following \textit{pseudo} negative log-likelihood:
\begin{eqnarray*}
\contrastL_{\indexset} (\K,d) & =& (1/2) \sum_{(j,k) \in\indexset}
\{ \dwt_{j,k}^2/(\K2^{2 d j} ) + \log(\K2^{2dj})\}\\
& =& \frac1{2 \K} \sum_{(j,k) \in\indexset}2^{-2 d j} \dwt_{j,k}^2 +
\frac{\cardinal{\indexset}}{2} \log\bigl(\K2^{2\jmean{\indexset}\,
d}\bigr) ,
\end{eqnarray*}
where $ \cardinal{\indexset}$ denotes the number of elements of the set
$\indexset$ and
$\jmean{\indexset}$ is defined as the average scale,
%
\begin{equation}\label{eq:etaDelta}
\jmean{\indexset}\eqdef\frac1{\cardinal{\indexset}}\sum_{(j,k)
\in\indexset} j .
\end{equation}
%
Define
$\hs^2_\indexset(d)\eqdef\argmin_{\K>0}\contrastL_{\indexset}
(\K,d) =
\cardinal{\indexset}^{-1} \sum_{(j,k) \in\indexset} 2^{-2 d j}
\dwt_{j,k}^2$.
The maximum pseudo-likelihood estimator of the memory parameter is
then equal to the minimum of the negative profile
log-likelihood (see~\cite{vandervaart1998}, page~403),
$\hd_\indexset\eqdef\break\argmin_{d\in\Rset} \contrastL_{\indexset}
(\hs^2_\indexset(d),d)$, that is,
%
\begin{equation} \label{eq:TildeJdef}
\hspace*{12pt}\hd_\indexset= \mathop\argmin_{d \in\Rset} \contrast_{\indexset
}(d),\qquad\mbox{where }
\contrast_{\indexset} (d) \eqdef\log\sum_{(j,k) \in\indexset}
2^{2d(\jmean{\indexset}-j)}\dwt_{j,k}^2 .
\end{equation}
%
If $\indexset$ contains at least two different scales, then $\contrast
_{\indexset} (d)\to\infty$ as
$d\to\pm\infty$ and thus
$\hd_\indexset$ is finite. The derivative of $\contrast_{\indexset}
(d)$ vanishes at
$d=\hd_{\indexset}$, that is, $\Sclt_\indexset(\hd_\indexset) = 0$, where
for all $d\in\Rset$,
%
\begin{equation}\label{eq:ScltDef}
\Sclt_\indexset(d)\eqdef\sum_{(j,k)\in\indexset} [j-\jmean
{\indexset}]
 2^{-2j d }\,\dwt_{j,k}^2 .
\end{equation}
We consider two specific choices for $\indexset$.
For any integers $n$, $j_0$ and $j_1$, $j_0 \leq j_1$, the set of all
available wavelet coefficients from $n$
observations $X_1,\dots,X_n$ having scale indices between $j_0$ and
$j_1$ is
%
\begin{equation}\label{eq:InDef}
\indexset_n(j_0,j_1) \eqdef\{ (j,k)\dvtx j_0 \leq j \leq j_1,
0\leq k < \nj{j} \} ,
\end{equation}
where $\nj{j}$ is given in \eqref{eq:Nj}.
Consider two sequences, $\{\nlowscale\}$ and $\{\nupscale\}$,
satisfying, for all~$n$,
%
\begin{equation}\label{eq:J0J1n}
0\leq\nlowscale< \nupscale\leq\nmaxscale, \qquad\nmaxscale\eqdef
\max
\{j\dvtx\nj{j}\geq1\} .
\end{equation}
The index $\nmaxscale$ is the maximal available scale index for the
sample size $n$;
$\nlowscale$ and $\nupscale$ will denote, respectively, the lower and
upper scale indices used in the pseudo-likelihood function.
The estimator will then be denoted $\hd_{\indexset_n(\nlowscale
,\nupscale)}$.
As shown \mbox{below}, in the semiparametric framework, the lower scale
$\nlowscale$ governs the rate of
convergence of $\hd_{\indexset_n(\nlowscale,\nupscale)}$ toward the
true memory parameter.
There are two possible settings as far as the upper scale $\nupscale$
is concerned:
\begin{longlist}[(S-1)]
\item[(S-1)]\label{it:caseJ1-J0fixed} $\nupscale-\nlowscale$ is fixed, equal
to $\k>0$;
\item[(S-2)]\label{it:caseJ1isJ} $\nupscale\leq\nmaxscale$ for all $n$ and
$\nupscale-\nlowscale\to\infty$ as $n\to\infty$.
\end{longlist}
\textup{(\hyperref[it:caseJ1-J0fixed]{S-1})} corresponds to using a fixed number of scales
and~\textup{(\hyperref[it:caseJ1isJ]{S-2})} corresponds to using a number
of scales tending to infinity.
We will establish the large sample properties of $\hd_{\indexset
_n(\nlowscale,\nupscale)}$ for these two cases.

The following theorem, proved in Section~\ref{sec:rates}, states that
under Condition~\ref{cond:GeneralCondition}, the estimator
$\hd_{\indexset_n(\nlowscale,\nupscale)}$ is consistent.
\begin{theorem}[(Rate of convergence)]
\label{thm:Rates}
Assume Condition~\textup{\ref{cond:GeneralCondition}}. Let $\{\nlowscale\}$ and
$\{\nupscale\}$ be two sequences
satisfying \textup{(\ref{eq:J0J1n})} and suppose that, as $n\to\infty$,
%
\begin{equation}
\label{eq:CondLn}
\nlowscale^2(n2^{-\nlowscale})^{-1/4}+\nlowscale^{-1}\to0  .
\end{equation}
The estimator $\hd_{\indexset_n(\nlowscale,\nupscale)}$ defined
by~\textup{(\ref{eq:TildeJdef})} and~\textup{(\ref{eq:InDef})} is then consistent with a
rate given by
%
\begin{equation}\label{eq:var-bias}
\hd_{\indexset_n(\nlowscale,\nupscale)}=\td+O_{\prob}\{
(n2^{-\nlowscale})^{-1/2}+ 2^{-\beta\nlowscale}\}  .
\end{equation}
\end{theorem}

By balancing the two terms in the bound~(\ref{eq:var-bias}), we obtain
the optimal rate.
\begin{corollary}[(Optimal rate)]
\label{coro:optimal-rate}
When $n\asymp2^{(1+2\beta)\nlowscale}$, we obtain the rate
%
\begin{equation}
\label{eq:optimalRate}
\hd_{\indexset_n(\nlowscale,\nupscale)}=\td+O_{\prob}
\bigl(n^{-\beta
/(1+2\beta)}\bigr) .
\end{equation}
\end{corollary}
\begin{pf*}{Proof}
By taking $n\asymp2^{(1+2\beta)\nlowscale}$, the condition
$\nlowscale
^{-1}+\nlowscale^2(n2^{-\nlowscale})^{-1/4}\to0$ is satisfied and
$(n L_n)^{-1/2} \asymp2^{-\beta L_n} \asymp n^{-\beta
/(1+2\beta)}$. This is the minimax rate
\cite{giraitisrobinsonsamarov1997}.
\end{pf*}
\begin{remark}
Observe that the setting of Theorem~\ref{thm:Rates} includes both
cases~\textup{(\hyperref[it:caseJ1-J0fixed]{S-1})} and~\textup{(\hyperref[it:caseJ1isJ]{S-2})}.
The difference between these settings will appear when computing the
limit variance in the Gaussian case; see
Theorem~\ref{thm:CLT} below.
\end{remark}

We shall now state a central limit theorem for the estimator $\hd
_{\indexset_n(\nlowscale,\nupscale)}$ of $\td$, under the additional
assumption that
$X$ is a Gaussian process. Extensions to non-Gaussian linear processes
will be considered in a future work.
We denote by $|\cdot|$ the Euclidean norm and define, for all $d\in
(1/2-\alpha,M]$ and $\dj\in\Nset$,
%
\begin{equation}\label{eq:bDint}
\intbdens[\psi]{\dj}{d} \eqdef
\int_{-\pi}^{\pi} |\bdensasymp[\psi]{\dj}{\lambda}{d}
|^2 \,\d\lambda
= (2\pi)^{-1}\sum_{\tau\in\Zset} \PCov^2\bigl(\dwt
_{0,0}^{(d)},\dwt
^{(d)}_{-\dj,\tau}\bigr) ,
\end{equation}
where we have used~(\ref{eq:CovFBMGenPsijk}). We denote, for all integer
$\k
\geq1$,
\begin{eqnarray}
\eta_\k&\eqdef&\sum_{j=0}^\k j \frac{2^{-j}}{2-2^{-\k}}
\quad\mbox{and}\quad
\kappa_\k\eqdef\sum_{j=0}^{\k} (j-\eta_{\k})^2 \frac
{2^{-j}}{2-2^{-\k}},\label{eq:eta-k}\\
\AsympVarWWE[\psi]{\td,\k} &\eqdef&\frac{\pi}{(2-2^{-\k})\kappa
_\k(\log(2)\Kvar[\psi]{\td})^2}\nonumber\\
&&{}\times\Biggl\{\intbdens[\psi]{0}{\td} +
\frac2{\kappa_\k}
\sum_{\dj=1}^\k\intbdens[\psi]{\dj}{\td}  2^{(2\td-1) \dj}\label{eq:varLimite-k}\\
&&\hspace*{86pt}{}\times\sum_{i=0}^{\k-\dj}\frac{2^{-i}}{2-2^{-\k}}(i-\eta_\k)(i+\dj-\eta
_\k)\Biggr\},\nonumber\\
\quad \AsympVarWWE[\psi]{\td,\infty} &\eqdef&
\frac{\pi}{[2\log(2)\Kvar[\psi]{\td}]^2} \Bigl\{\intbdens[\psi
]{0}{\td} +
2\sum_{\dj=1}^\infty\intbdens[\psi]{\dj}{\td}  2^{(2\td-1)
\dj}\Bigr\} ,\label{eq:varLimite}
\end{eqnarray}
where $\Kvar[\psi]{d}$ is defined in \eqref{eq:Kpsi}. The following
theorem is proved in Section~\ref{sec:rates}.
\begin{theorem}[(CLT)]\label{thm:CLT}
Let $X$ be a Gaussian process having generalized spectral measure~\textup{(\ref{eq:fmodele})}
with $\td\in\Rset$ and $\mes^\ast\in\calH(\beta,\gamma
,\varepsilon)$
with $\mes^\ast(-\varepsilon,\varepsilon)>0$, where
$\gamma>0$, $\beta\in(0,2]$ and $\varepsilon\in(0,\pi]$.
Let $\{\nlowscale\}$ be a sequence such that
%
\begin{equation}\label{eq:CondBiaisNegligeable}
\nlowscale^2(n2^{-\nlowscale})^{-1/4}+ n 2^{- (1+2\beta) \nlowscale}
\to0
\end{equation}
and $\{\nupscale\}$ be a sequence such that either~\textup{(\hyperref[it:caseJ1-J0fixed]{S-1})}
or~\textup{(\hyperref[it:caseJ1isJ]{S-2})} holds.
Then, under \allWA, we have, as $n\to\infty$,
%
\begin{equation}\label{eq:clt}
(n2^{-\nlowscale})^{1/2}\bigl(\hd_{\indexset_n(\nlowscale,\nupscale
)}-\td
\bigr)\cl\calN[ 0,\AsympVarWWE[\psi]{\td,\k} ]  ,
\end{equation}
where $\k=\lim_{n\to\infty}(\nupscale-\nlowscale)\in\{1,2,\dots
,\infty\}$.
\end{theorem}
\begin{remark}
The condition~(\ref{eq:CondBiaisNegligeable}) is similar to~(\ref
{eq:CondLn}), but ensures, in addition, that the bias
in~(\ref{eq:var-bias}) is asymptotically negligible.
\end{remark}
\begin{remark}\label{rem:W5}
The larger the value of $\beta$, the smaller the size of the allowed
range for $\td$ in~(\hyperref[item:Wvstd]{W-5}) for a given
decay exponent $\alpha$ and number $M$ of vanishing moments. Indeed, the
range in~(\hyperref[item:Wvstd]{W-5})
has been chosen so as to obtain a bound on the bias which corresponds
to the best possible rate under the condition
$\mes^\ast\in\calH(\beta,\gamma,\varepsilon)$. If~(\hyperref[item:Wvstd]{W-5}) is
replaced by the weakest condition $\td\in(1/2-\alpha,M]$,
which does not depend on $\beta$, the same CLT~(\ref{eq:clt}) holds, but
$\beta$ in condition~(\ref{eq:CondBiaisNegligeable})
must be replaced by $\beta'\in(0,\beta]$. This $\beta'$ must satisfy
$1/2-\alpha<(1+\beta')/2-\alpha<\td$, that is, $0<\beta'<2(\td
+\alpha
)-1$. When $\beta'<\beta$, one gets a slower rate in~(\ref{eq:clt}).
\end{remark}
\begin{remark}\label{rem:compare}
Relation~(\ref{eq:clt}) holds under~\textup{(\hyperref[it:caseJ1-J0fixed]{S-1})}, where $\k
<\infty$ and~\textup{(\hyperref[it:caseJ1isJ]{S-2})}, where $\k=\infty$.
It follows from~(\ref{eq:iuBound}) and~(\ref{eq:Varkinfty}) that
$\AsympVarWWE[\psi]{\td,\k}
\to\AsympVarWWE[\psi]{\td,\infty}<\infty$ as $\k\to\infty$. Our
numerical experiments suggest that in some cases, one may have
$\AsympVarWWE[\psi]{\td,\k}\leq\AsympVarWWE[\psi]{\td,\k'}$
with $\k
\leq\k'$; see the bottom left panel of
Figure~\ref{fig:numvar}. In that figure, one indeed notices a bending
of the curves for large $d$, which is more pronounced
for small values of $M$ and may be due to a correlation between the
wavelet coefficients across scales.
\end{remark}
\begin{remark}\label{rem:J1-J0fini}
The most natural choice is $\nupscale=\nmaxscale$, which amounts to
using
all the available wavelet
coefficients with scale index larger than $\nlowscale$.
The case~\textup{(\hyperref[it:caseJ1-J0fixed]{S-1})} is nevertheless of interest.
In practice, the number of observations $n$ is finite and the number of
available scales $\nmaxscale-\nlowscale$ can be small.
Since, when $n$ is finite, it is always possible to interpret the
estimator $\hd_{\indexset_n(\nlowscale,\nmaxscale)}$ as $\hd
_{\indexset
_n(\nlowscale,\nlowscale+\k)}$
with $\k=\nmaxscale-\nlowscale$, one may approximate the
distribution of
$(n2^{-\nlowscale})^{1/2}(\hd_{\indexset_n(\nlowscale,\nmaxscale
)}-\td)$ either by
$\calN(0,\AsympVarWWE[\psi]{\td,\k})$ or by
$\calN(0,\AsympVarWWE[\psi]{\td,\infty})$. Since the former involves
only a single limit, it is likely to provide a better
approximation for finite $n$. Another interesting application
involves considering online estimators of $\td$: online
computation of wavelet coefficients is easier when the number of scales
is fixed; see~\cite{roughanveitchabry2000}.
\end{remark}

\section{Discussion}
\label{sec:disc-concl-remarks}
The asymptotic variance $\AsympVarWWE[\psi]{d,\k}$ is defined for all
$\k\in\{1,2,\break\dots,\infty\}$ and all
$1/2+\alpha<d\leq M$ by~(\ref{eq:varLimite-k}) and~(\ref{eq:varLimite}).
Its expression involves the range of scales $\k$ and
the $L^2$-norm $\intbdens[\psi]{\dj}{\td}$ of the asymptotic spectral
density $\bdensasymp[\psi]{\dj}{\lambda}{d}$ of the
wavelet coefficients, both for the ``within'' scales ($\dj=0$) and the
``between'' scales ($\dj>0$). The choice of wavelets does not
matter much, as Figure~\ref{fig:numvar} indicates. One can use
Daubechies wavelet or Coiflets (for which the scale function
also has vanishing moments). What matters is the number of vanishing
moments $M$ and the decay exponent $\alpha$, which both
determine the frequency resolution of $\psi$. For wavelets derived from
a multiresolution analysis, $M$ is always known and
\cite{cohen2003}, Remark~2.7.1, page~86, provides a sequence of lower
bounds tending to $\alpha$ (we used such lower bounds
for the Coiflets used below). For the Daubechies wavelet with $M$
vanishing moments, an analytic formula giving $\alpha$
is available; see \cite{daubechies1992}, equation~(7.1.23), page~225 and the table on
page~226, and
note that our $\alpha$ equals the $\alpha$ of \cite{daubechies1992}
plus~1.

\subsection{The ideal Shannon wavelet case}
\begin{figure}

\includegraphics{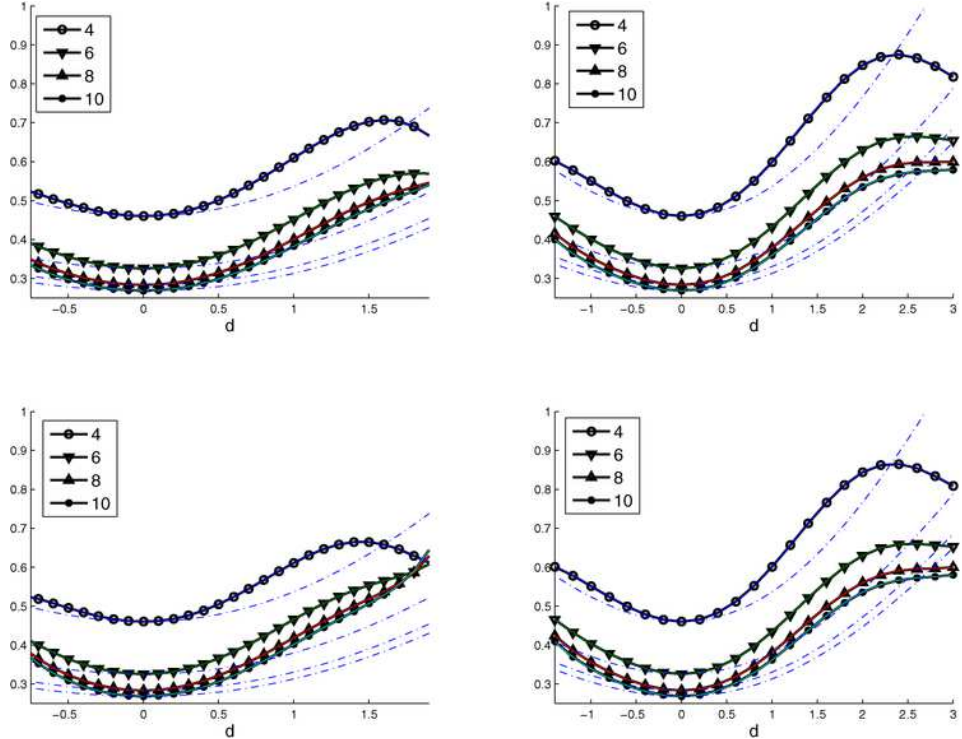}

\caption{Numerical computations of the asymptotic variance
$\AsympVarWWE[\psi]{d,\k}$ for the Coiflets and
Daubechies wavelets for different values of the number of scales $\k
=4, 6, 8, 10$ and of the number of vanishing moments
$M=2,4$. Top row: Coiflets; bottom row: Daubechies wavelets;
left column: $M=2$; right column: $M=4$. The dash-dot lines are the
asymptotic variances for the Shannon wavelet
[see~(\protect\ref{eq:shannonApprox})] with $\k=4,6,8,10$. For a given $\k$,
the variances for different orthogonal wavelets
coincide at $d=0$; see the comment following~(\protect\ref{eq:univLB}). The
right and left columns have different horizontal
scales because different values of $M$ yield different ranges for $d$.}
\label{fig:numvar}
\end{figure}

The so-called Shannon wavelet $\psi_S$ is such that its Fourier
transform $\hat{\psi}_S$ satisfies $|\hat{\psi}_S(\xi)|^2=1$
for $|\xi|\in[\pi,2\pi]$ and is zero otherwise. This wavelet
satisfies~(\hyperref[item:psiHat]{W-2})--(\hyperref[item:MIM]{W-4}) for arbitrary large $M$
and~$\alpha$, but
does not have compact support, hence it does not satisfy~(\hyperref[item:Wreg]{W-1}).
We may not, therefore, choose this wavelet in our analysis.
It is of interest, however, because it gives a rough idea of what
happens when $\alpha$ and $M$ are large since one can always
construct a wavelet $\psi$ satisfying~(\hyperref[item:Wreg]{W-1})--(\hyperref[item:MIM]{W-4})
which is arbitrarily close to the Shannon wavelet.
Using the Shannon wavelet in~(\ref{eq:bDpsi}), we get, for all
$\lambda
\in(-\pi,\pi)$,
$\bdensasymp[\psi_S]{\dj}{\lambda}{d}=0$ for $\dj\geq1$ and
$\bdensasymp[\psi_S]{0}{\lambda}{d}=(2\pi-|\lambda|)^{-2d}$ so that,
for all $d\in\Rset$,~(\ref{eq:varLimite-k}) becomes
%
\begin{equation}
\label{eq:shannonApprox}
\hspace*{16pt}\AsympVarWWE[\psi_S]{d,\k} = \frac{\pi g(-4d)}
{2(2-2^{-\k})\kappa_\k\log^2(2)  g^2(-2d)}
\qquad\mbox{where } g(x)=\int_\pi^{2\pi} \lambda^x \,\d
\lambda.\hspace*{-10pt}
\end{equation}
This $\AsympVarWWE[\psi_S]{d,\k}$ is displayed in Figure~\ref{fig:numvar}.

\subsection{Universal lower bound for $\mathrm{I}_{0}(d)$}

For $\k=\infty$, using the facts that $\intbdens[\psi]{0}{d}\geq0$ for $\dj
\geq
1$ and,
by the Jensen inequality in~(\ref{eq:bDint}), $\intbdens[\psi
]{0}{d}\geq\Kvarsq[\psi]{d}/(2\pi)$, we have, for all
$1/2+\alpha<d\leq M$,
%
\begin{equation}
\label{eq:univLB}
\AsympVarWWE[\psi]{d,\infty} \geq(8\log^2(2)
)^{-1}\simeq
0.2602  .
\end{equation}
This inequality is sharp when $d=0$ and the wavelet family $\{\psi
_{j,k}\}_{j,k}$ forms an orthonormal basis. This is because, in this
case, the lower bound
$(8 \log^2(2))^{-1}$ in \eqref{eq:univLB} equals $\AsympVarWWE[\psi
]{0,\infty}$. Indeed, by~(\ref{eq:CovFBMGen}) and Parsevals theorem,
the wavelet coefficients $\{B_{(0)}(\psi_{j,k})\}_{j,k}$ are a centered
white noise with variance $2\pi$ and, by~(\ref{eq:Kpsi}) and~(\ref
{eq:bDint}),
$\Kvar[\psi]{0}=2\pi$ and $\intbdens[\psi]{\dj}{0}=2\pi\1(\dj
=0)$. Then,
$\AsympVarWWE[\psi]{0,\k} = (2(2-2^{-\k})\kappa_\k\log
^2(2))^{-1}$.
Since $\kappa_\k$ is increasing with $\k$ and tends to $2$ as $\k
\to
\infty$ (see Lemma \ref{lem:etavareta}),
$\AsympVarWWE[\psi]{0,\k}\geq(8\log^2(2))^{-1} =
\AsympVarWWE[\psi]{0,\infty}$.
Hence, the lower bound~(\ref{eq:univLB}) is attained at $\td=0$ if $\{
\psi
_{j,k}\}_{j,k}$ is an orthonormal
basis.


\subsection{Numerical computations}
For a given wavelet $\psi$, we can compute the variances
$\AsympVarWWE[\psi]{d,\k}$ numerically for any $\k=1,2,\dots
,\infty$
and $1/2+\alpha<d\leq M$. It is easily shown that
$d\mapsto\AsympVarWWE[\psi]{d,\k}$ is infinitely differentiable on
$1/2+\alpha<d\leq M$ so that interpolation can be
used between two different values of $d$. We compared numerical values
of $\AsympVarWWE[\psi]{d,\k}$ for four different
wavelets, with $\k=4,6,8,10$, and compared them with the Shannon
approximation~(\ref{eq:shannonApprox}); see
Figure~\ref{fig:numvar}. We used as wavelets two Daubechies wavelets
which have $M=2$
and $M=4$ vanishing moments, and $\alpha=1.3390$ and $\alpha=1.9125$
decay exponents, respectively, and two so-called Coiflets
with the same number of vanishing moments, and $\alpha>1.6196$ and
$\alpha>1.9834$ decay exponents respectively.
For a given number $M$ of vanishing moments, the Coiflet has a larger
support than the Daubechies wavelet, resulting in a better
decay exponent. The asymptotic variances are different for $M=2$, in
particular, for negative $d$'s, the
Coiflet asymptotic variance is closer to that of the Shannon wavelet. The
asymptotic variances are very close for $M=4$.

\subsection{Comparison with Fourier estimators}
Semiparametric Fourier estimators are based on the periodogram.
To allow comparison with Fourier estimators, we must first link the
normalization factor $n2^{-\lowscale_n}$ with the bandwidth
parameter $m_n$ (the index of the largest normalized frequency) used by
semiparametric Fourier estimators.
A Fourier estimator with bandwidth $m_n$ projects the observations
$[X_1\dots X_n]^T$ on the space generated by the vectors
$\{\cos(2\pi k \cdot/n),\sin(2\pi k \cdot/n) \}$, $k=1,\dots,m_n$, whose
dimension is $2m_n$; on the other hand, the wavelet coefficients $\{
\dwt
_{j,k}, j\geq\lowscale, k=0,\dots,\nj{j}-1\}$ used in the
wavelet estimator correspond to a projection on a space whose dimension
is at most
$\sum_{j=\nlowscale}^{\nmaxscale}\nj{j}\sim2n2^{-\lowscale_n}$, where
the equivalence holds as $n\to\infty$ and $n2^{-\lowscale_n}\to
\infty$,
by applying~(\ref{eq:etavaretafirst}) with $j_0=\nlowscale$,
$j_1=\nmaxscale$ and $p=1$.
Hence, for $m_n$ or
$n2^{-\lowscale_n}$ large, it makes sense to consider $n2^{-\lowscale
_n}$ as an analog of the bandwidth parameter $m_n$.
The maximal scale index $\upscale_n$ is similarly related to the \textit
{trimming number} (the index of the smallest normalized frequency), often
denoted by $l_n$ (see~\cite{robinson1995l}), that is, $l_n \sim n
2^{-\upscale_n}$. We stress that, in absence of trends, there is no
need to trim coarsest scales.

With the above notation, the assumption~(\ref{eq:CondLn}) in
Theorem~\ref{thm:Rates} becomes
$m_n/n+(\log n/m_n)^8 m_n^{-1}\to0$ and the conclusion~(\ref
{eq:var-bias}) is expressed as
$\hd=\td+\break O_\prob(m_n^{-1/2}+(m_n/n)^\beta)$. The
assumption~(\ref{eq:CondBiaisNegligeable}) becomes
$(\log n/m_n)^8 m_n^{-1}+m_n^{1+2\beta}/n^{2\beta}\to0$ and the rate
of convergence in~(\ref{eq:clt}) is $m_n^{1/2}$.

The most efficient Fourier estimator is the local Whittle (Fourier)
estimator studied in~\cite{robinson1995g}; provided
that
\begin{longlist}[(1)]
\item[(1)]\label{item:GSE:condition-dsp} the process $\{ X_k \}$ is
stationary and has spectral $f(\lambda)= |1-\rme^{-\rmi
\lambda}|^{-2 \td} f^\star(\lambda)$ with $\td\in(-1/2,1/2)$ and
$f^\ast(\lambda)=f^\ast(0)+O(|\lambda|^\beta)$ as $\lambda\to0$,
\item[(2)]\label{item:GSE:condition-linear} the process $\{ X_k \}$ is
linear and causal, $X_k = \sum_{j=0}^\infty a_j Z_{k-j}$,
where $\{ Z_k \}$ is a martingale increment sequence satisfying $\PE
[Z_k^2 \mid\mathcal{F}_{k-1}]= 1$ a.s.,
$\PE[Z_k^3 | \mathcal{F}_{k-1}]= \mu^3$ a.s. and $\PE[Z_k^4]=\PE
[Z_1^4] $, where $\mathcal{F}_k= \sigma(Z_{k-l},l\geq0)$
and $a(\lambda) \eqdef\sum_{k=0}^\infty a_k \rme^{- \rmi k \lambda}$
is differentiable in a neighborhood $(0,\delta)$ of
the origin and
$|\d a / \d\lambda(\lambda)| =\break O( |a(\lambda)|/ \lambda)$
as $\lambda\to0^+$ (see~A2$'$)
\item[(3)]\label{item:GSE:condition-bandwidth} $m_n^{-1} + (\log m_n)^2
m_n^{1+2\beta} /n^{2 \beta} \to0$ (see A4$'$) ,
\end{longlist}
then $m_n^{1/2} (\hat{d}_{m_n} - \td)$ is asymptotically zero-mean
Gaussian with variance 1$/$4.
This asymptotic variance is smaller than (but very close to) our lower
bound in~(\ref{eq:univLB}) and comparable to the asymptotic variance obtained
numerically for the Daubechies wavelet with two vanishing moments; see
the left-hand panel in Figure~\ref{fig:numvar}.
Also, note that while the asymptotic variance of the Fourier estimators
is a constant, the asymptotic variances of the wavelet estimators
depend on $\td$ (see Figure~\ref{fig:numvar}). In practice, one
estimates the limiting variance $\AsympVarWWE[\psi]{\td,\k}$
by $\AsympVarWWE[\psi]{\hd,\k}$ in order to construct asymptotic
confidence intervals. The continuity of
$\AsympVarWWE[\psi]{\cdot,\k}$ and the consistency of $\hd$ justify
this procedure.

We would like to stress, however, that the wavelet estimator has some
distinctive advantages.
From a theoretical standpoint, for a given $\beta$, the wavelet
estimator is rate optimal, that is, for $\beta\in(0,2]$,
the rate is $n^{\beta/1+2\beta}$ (see Corollary~\ref{coro:optimal-rate}) and the CLT is obtained for any rate
$o(n^{\beta/1+2\beta})$. For the local Whittle Fourier estimator, the
best rate of convergence is $O((n/\log^2(n))^{\beta/1+2\beta})$ and the
CLT is obtained for any rate $o((n/\log^2(n))^{\beta/1+2\beta})$. This
means that for any given $\beta$, the wavelet estimator has a faster
rate of convergence and can therefore yield, for an appropriate
admissible choice of the finest scale, shorter confidence intervals.
Another advantage of the wavelet Whittle estimator over this estimator
is that the optimal rate of convergence is shown to
hold for $\mes^\star\in\calH(\beta,\gamma,\varepsilon)$ without any
further regularity assumption, such as the
density $f^\ast$ of $\mes^\ast$ having to be differentiable in a neighborhood
of zero, with a given growth of the logarithmic
derivative.
To the best of our knowledge, the GPH estimator is the only
Fourier estimator which has been shown, in a Gaussian context, to
achieve the rate $O(n^{\beta/(1+2\beta)})$
(see~\cite{giraitisrobinsonsamarov1997}); its asymptotic variance is
$\pi^2/24\simeq0.4112$. It is larger than the lower bound~(\ref{eq:univLB}) and larger than
the asymptotic variance obtained by using standard Daubechies wavelets
with $\k\geq6$ on the
range $(-1/2,1/2)$ of $\td$ allowed for the GPH estimator (see
Figure~\ref{fig:numvar}). When pooling frequencies, the
asymptotic variance of the GPH estimator improves and tends to $1/4$
(the local Whittle Fourier asymptotic variance) as the
number of pooled frequencies tends to infinity; see~\cite{robinson1995l}.
%
%

Thus far, we have compared our local Whittle wavelet estimator with the
local Whittle Fourier (\LWF) and GPH estimators in the context of a stationary
and invertible process $X$, that is, for $\td\in(-1/2,1/2)$.
As already mentioned, the wavelet estimators can be used for
arbitrarily large ranges of the parameter $\td$ by appropriately choosing
the wavelet so that (\hyperref[item:Wvstd]{W-5}) holds.
There are two main ways of adapting the \LWF\ estimator to larger
ranges of $d$: differentiating and
tapering the data (see~\cite{velasco1999}) or, as promoted by
\cite{shimotsuphillips2005}, modifying the
local Whittle likelihood, yielding the so-called exact local Whittle
Fourier (ELWF) estimator.
The theoretical analysis of these methods is performed under the same
set of assumptions as in \cite{robinson1995g}, so the same
comments on the
nonoptimality of the rate and on the restriction on $f^\star$ apply.
Also, note that the model considered by~\cite{shimotsuphillips2005}
for $X$ differs from the model of integrated processes
defined by \eqref{eq:linearProc} and is not time-shift invariant; see
their equation~(1). In addition, their estimator is not invariant under
the addition of a constant in the data, a drawback which is not easily
dealt with; see their Remark~2.
The asymptotic variance of the \ELWF\ estimator has been shown to be
$1/4$, the same as the LFW estimator, provided that the range $(\Delta
_1,\Delta_2)$ for $\td$ is of width
$\Delta_2-\Delta_1 \leq9/2$. The asymptotic variance of our local
Whittle wavelet estimator with eight scales, using the
Daubechies wavelet with $M=4$ zero moments, is at most $0.6$ on a range
of same width;
see the left-hand panel in Figure~\ref{fig:numvar}. Again, this
comparison does not take into account the logarithmic factor in the
rate of convergence imposed by the conditions
on the bandwidth $m_n$. Concerning the asymptotic variances of tapered
Fourier estimators,
increasing the allowed range for $\td$ means increasing the taper order
(see \cite{hurvichmoulinessoulier2002} and
\cite{robinsonhenry2003}), which, as already explained, inflates the
asymptotic variance of the estimates.
In contrast, for the wavelet methods, by increasing the number of
vanishing moments $M$ of, say, a Daubechies wavelet, the allowed range
for $\td$ is arbitrarily large
while the asymptotic variance converges to the ideal Shannon wavelet
case, derived in \eqref{eq:shannonApprox}; the numerical values are
displayed in Figure~\ref{fig:numvar}
for different values of the number of scales $\k$. The figure shows
that larger values of $\k$ tend to yield a smaller
asymptotic variance. One should thus choose the largest possible $M$
and the maximal number of scales. This prescription
cannot be applied to a small sample because increasing
the support of the wavelet decreases the number of available scales.
The Daubechies wavelets with $M=2$ to $M=4$ are commonly used in practice.

%
%
From a practical standpoint, the wavelet estimator is computationally
more efficient than the aforementioned Fourier estimators.
Using the fast pyramidal algorithm, the wavelet transform coefficients
are computed in $O(n)$ operations. The function $d
\mapsto\contrast_{\indexset}(d)$ can be minimized using the Newton
algorithm \cite{boydvandenberghe2004}, Chapter~9.5,
whose convergence is guaranteed because $\contrast_{\indexset}(d)$ is
convex in $d$.
The complexity of the minimization procedure is related to the
computational cost of evaluation of the function $\contrast_{\indexset
}$ and its two first derivatives.
Assume that these functions need to be evaluated at $p$ distinct values
$d_1,\dots,d_p$.
We first compute the empirical variance of the wavelet coefficients
$\nj
{j}^{-1}\sum_{k=0}^{\nj{j}-1}\dwt_{j,k}^2$ for the scales $j\in\{
\lowscale_n, \dots,\upscale_n\}$, which does not depend
on $d$ and requires $O(n)$ operations. For $\indexset= \indexset
_n(\lowscale_n,\upscale_n)$,
$\contrast_{\indexset}$ and all of its derivatives are linear
combinations of these $\upscale_n-\lowscale_n+1= O(\log(n))$ empirical
variances with weights depending on $d$.
The total complexity for computing the wavelet Whittle estimator in an
algorithm involving
$p$ iterations is thus $O(n + p \log(n))$.
The local Whittle Fourier (\LWF) contrast being convex, the same Newton
algorithm converges, but the complexity is slightly
higher. The computation of the Fourier coefficients requires $O(n \log
(\bolds{n}))$ operations.
The number of terms in the \LWF\ contrast function (see \cite{robinson1995g},
page~1633)
is of order $m_n$ [which is typically of order
$O(n^{\gamma})$, where $\gamma\in(0, 1/1+2\beta)$],
so the evaluation of the \LWF\ contrast function (and its derivatives)
for $p$ distinct values of the memory parameter $d_1, \dots, d_p$
requires $O(p m_n)$ operations.
The overall complexity of computing the \LWF\ estimator in a Newton algorithm involving $p$
steps is therefore $O( n \log(n) + p m_n)$.
Differentiating and tapering the data only adds $O(n)$ operations, so
 the same complexity applies in this case.
The \ELWF\ estimator is much more computationally demanding and is
impractical for large data sets:
for each value of the memory coefficient $d$ at which the
pseudo-likelihood function is evaluated, the algorithm calls for
the fractional integration or differentiation of the observations,
namely, $(\Delta^d X)_k,k=1,\dots,n$, and the computation of the Fourier
transform of $\{ (\Delta^d X)_1, \dots, (\Delta^d X)_n \}$. In this
context, $(\Delta^d X)_k\eqdef\sum_{l=0}^k \frac{(-d)_l}{l!} X_{k-l}$,
$k=1,\dots,n$, where $(x)_0=1$ and $(x)_k=x(x+1)\cdots(x+k-1)$ for
$k\geq
1$ denote the
Pochhammer symbols.
The complexity of this procedure is thus $O(n^2 + n \log(n))$. The
complexity for $p$ function evaluations, therefore,
is $O(p (n^2 + n \log(n)))$. The convexity of the criterion is not
assured, so a minimization algorithm can possibly be trapped in a
local minimum.
These drawbacks make the \ELWF\ estimator impractical for large data
sets, say of size $10^6-10^7$, as encountered in teletraffic analysis or
high-frequency financial data.

\section{\texorpdfstring{Condition~\protect\ref{cond:GeneralCondition} holds for linear and
Gaussian processes}{Condition~1 holds for linear and
Gaussian processes}}\label{sec:assumpcons}
\mbox{}

\begin{pf*}{Proof of Theorem~\protect\ref{thm:AssumpCons}}
For any scale index $j \in\Nset$, define by $\{ h_{j,l} \}_{l \in
\Zset
}$ the sequence $h_{j,l} \eqdef2^{-j/2}
\int_{-\infty}^\infty\phi(t+l)\psi(2^{-j}t)\,\d t$ and by
$H_j(\lambda) \eqdef\sum_{l\in\Zset}h_{j,l}\rme^{-\rmi\lambda
l}$ its
associated discrete-time Fourier transform.
Since $\phi$ and $\psi$ are compactly supported, $\{ h_{j,l} \}$ has a
finite number of nonzero coefficients.
As shown by~\cite{moulinesrouefftaqqu2007a}, Relation~13, for any
sequence $\{ x_l \}_{l \in\Zset}$,
the discrete wavelet transform coefficients at scale $j$ are given by
$\dwt^x_{j,k}=\sum_{l\in\Zset} x_lh_{j,2^jk-l}$. In
addition, it follows from ~\cite{moulinesrouefftaqqu2007a},
Relation~16, that
$H_j(\lambda)= (1 - \rme^{-\rmi\lambda})^M \tilde{H}_j(\lambda)$,
where $\tilde{H}_j(\lambda)$ is a trigonometric polynomial, that is,
$\tilde
{H}_j(\lambda) = \sum_{l \in\Zset} \tilde{h}_{j,l} \rme^{-\rmi
\lambda l}$, where $\{ \tilde{h}_{j,l} \}$ has
a finite number of nonzero coefficients.

Define $\overline{\mes}$ and $\tilde{\mes}$ as the restrictions of
$\mes
$ on $[-\varepsilon,\varepsilon]$ and
on its complementary set, respectively. These definitions imply that
%
\begin{equation}
\label{eq:decompositionvj}
\sigma^{2}_{j}(\nu)=\sigma^{2}_{j}(\overline{\nu})+\sigma^{2}_{j}(\tilde{\nu}).
\end{equation}
Since $\mes^\ast\in\calH(\beta,\gamma,\varepsilon)$, the corresponding
decomposition for $\mes^\ast$ reads as in~(\ref{eq:decompmesstar}),
so $\overline{\mes}$ admits a density
$f(\lambda)=|1-\rme^{-\rmi\lambda}|^{-2\td}f^\ast(\lambda)$ on
$\lambda
\in[-\pi,\pi]$, where $f^\ast(\lambda)=0$ for
$\lambda\notin[-\varepsilon,\varepsilon]$ and $|f^\ast(\lambda
)-f^\ast
(0)|\leq\gamma f^\ast(0)|\lambda|^\beta$ on
$\lambda\in[-\varepsilon,\varepsilon]$. Hence, (\ref{eq:DjApprox})
holds: by \cite{moulinesrouefftaqqu2007a}, Theorem~1,
there exists a constant $C$ such that for all $j \geq0$ and $\lambda
\in(-\pi,\pi)$,
%
\begin{equation}\label{eq:DjApproxbis}
|\bdens[\phi,\psi]{j,\dj}{\lambda}{f} - f^\ast(0)
\bdensasymp
[\psi]{\dj}{\lambda}{\td}  2^{2j \td}|
\leq C f^\ast(0)  \bar{\gamma}  2^{(2\td-\beta)j}.
\end{equation}
Recall that $\bdens[\phi,\psi]{j,0}{\lambda}{f}$ is the spectral
density of a stationary series with variance $
\sigma^{2}_{j}(\overline{\nu})=\int_{-\pi}^\pi\bdens[\phi,\psi
]{j,0}{\lambda
}{f}\,\d\lambda$.
Similarly, by~(\ref{eq:CovFBMGenPsijk}) and~(\ref{eq:Kpsi}),
$\bdensasymp[\psi]{0}{\lambda}{\td}$ is the spectral density
of a stationary series with variance $\Kvar[\psi]{\td}$. Thus, after
integration on $\lambda\in(-\pi,\pi)$,~(\ref{eq:DjApproxbis}) with
$\dj=0$
yields
%
\begin{equation}
\label{eq:DjApproxbis1}
| \sigma^{2}_{j}(\overline{\nu}) - f^\ast(0)  \Kvar[\psi]{\td
}
2^{2j\td} | \leq2\pi C f^\ast(0)
\bar{\gamma}  2^{(2\td-\beta)j} .
\end{equation}
By \cite{moulinesrouefftaqqu2007a}, Proposition 9, there exists a
constant $C$ such that $|H_j(\lambda)| \leq C
2^{j(M+1/2)}\times\break |\lambda|^M (1+2^j |\lambda|)^{-\alpha-M}$ for any
$\lambda\in[-\pi,+\pi]$, which implies that
\begin{eqnarray}
\label{eq:vjtildefBound}
\hspace*{24pt}\sigma^{2}_{j}(\tilde{\nu}) &=&2 \int_{\varepsilon}^\pi|H_j(\lambda
)|^2
\mes(\d\lambda)
\leq C2^{(1+2M)j} \int_{\varepsilon}^\pi\lambda^{2M}(1+2^j\lambda
)^{-2\alpha-2M}  \mes(\d\lambda)\nonumber\\
&\leq& C\pi^{2M}  2^{(1+2M)j}(1+\varepsilon2^j)^{-2\alpha-2M}
\mes([\varepsilon,\pi])\\
&=& O\bigl(2^{j(1-2\alpha)}\bigr)
=o\bigl(2^{j(2\td-\beta)}\bigr) ,\nonumber
\end{eqnarray}
since, by~(\hyperref[item:Wvstd]{W-5}), $1-2\alpha-2\td+\beta<0$.
Relations~\eqref{eq:decompositionvj},
\eqref{eq:DjApproxbis1} and \eqref{eq:vjtildefBound} prove (\ref{eq:BiasTermAssump}).

We now consider~(\ref{eq:rosenthalAssump}). We have, for all $j\geq0$
and $n \geq1$, (see \cite{rosenblatt1985}, Theorem~2, page~34),
\begin{eqnarray}
\label{eq:emvarLinCase}
\PVar\Biggl(\sum_{k=0}^{\nj{j}-1} \dwt_{j,k}^{2}\Biggr)&=&\sum
_{\tau=-\nj
{j}+1}^{\nj{j}-1}(\nj{j}-|\tau|)
\PCov(\dwt_{j,0}^2,\dwt_{j,\tau}^2)\nonumber\\
&=&\sum_{\tau=-\nj{j}+1}^{\nj{j}-1} (\nj{j}-|\tau|)
[2\PCov^2(\dwt_{j,0},\dwt_{j,\tau})\\
&&{}+\Cum(\dwt_{j,0},\dwt
_{j,0},\dwt
_{j,\tau},\dwt_{j,\tau})].\nonumber
\end{eqnarray}
Using~(\ref{eq:linearProc}), since $M \geq k_0$, we may write
%
\begin{equation}
\label{eq:bjk}
\dwt_{j,k}=\sum_{t\in\Zset} \tilde{h}_{j,2^jk-t} (\diffop^M X)_t
= \sum_{t\in\Zset} b_{j,2^jk-t} Z_t,
\end{equation}
where $b_{j,\bolds{\cdot}}\eqdef\tilde{h}_{j,\bolds{\cdot}}\star(\diffop^{M-k_0} a)$
belongs to $\ell^2(\Zset)$. By~(\ref{eq:cumFour}), we thus obtain
\[
\Cum(\dwt_{j,0},\dwt_{j,0},\dwt_{j,\tau},\dwt_{j,\tau})=(\PE
[Z_1^4]-3)\sum_{t\in\Zset} b_{j,t}^2b_{j,2^j\tau-t}^2,
\]
which, in turns, implies that
%
\begin{eqnarray}
\label{eq:cumfourBound}
\sum_{\tau\in\Zset}|\Cum(\dwt_{j,0},\dwt_{j,0},\dwt
_{j,\tau},\dwt
_{j,\tau})|
&=&|\PE[Z_1^4]-3|\sum_{t,\tau\in\Zset} b_{j,t}^2b_{j,2^j\tau
-t}^2\nonumber\\[-8pt]\\[-8pt]
&\leq&|\PE[Z_1^4]-3|\sigma^{4}_{j}(\nu)\nonumber
\end{eqnarray}
since, by~(\ref{eq:bjk}), $\sum_{t}b_{j,t}^2=\sigma^{2}_{j}(\nu)$.

We shall now bound $\sum_{\tau=-\nj{j}+1}^{\nj{j}-1}
\PCov^2(\dwt_{j,0},\dwt_{j,\tau})$.
One can define uncorrelated wavelet coefficients $\{ \overline{\dwt
}_{j,k} \}$ and $\{ \tilde{\dwt}_{j,k} \}$, associated
with the generalized spectral measures $\overline{\mes}$ and $\tilde
{\mes}$, respectively and such that $\dwt_{j,k}=
\overline{\dwt}_{j,k}+\tilde{\dwt}_{j,k}$ for all $j\geq0$ and
$k\in
\Zset$.
Therefore,
$\PCov^2(\dwt_{j,0},\dwt_{j,\tau})= \PCov^2(\overline{\dwt
}_{j,0},\overline{\dwt}_{j,0}) + \PCov^2(\tilde{\dwt}_{j,0},\tilde
{\dwt
}_{j,\tau})
+2 \PCov(\overline{\dwt}_{j,0},\overline{\dwt}_{j,\tau}) \PCov
(\tilde
{\dwt}_{j,0},\tilde{\dwt}_{j,\tau})$.
By~(\ref{eq:BiasTermAssump}), $\sigma^{2}_{j}(\nu)\asymp2^{2j\td}$.
Therefore, by~(\ref{eq:DjApproxbis}) and using
\cite{moulinesrouefftaqqu2007a}, Proposition~3, equation~(30), for all $j\geq0$,
$\{\sigma^{-1}_{j}(\nu) \overline{\dwt}_{j,k} ,k\in\Zset\}$ is a
stationary process whose spectral density is bounded
above by a constant independent of $j$. Parsevals theorem implies
that $\sup_{j\geq1} \sigma^{-4}_{j}(\nu)
\sum_{\tau\in\Zset}\PCov^2(\overline{\dwt}_{j,0},\overline{\dwt
}_{j,\tau
})<\infty$, hence
%
\begin{equation}
\label{eq:cuvSqrBound1}
\sup_{n \geq1} \sup_{j=1,\dots,J_n} \nj{j}^{-1} \sigma^{-4}_{j}(\nu)
\sum
_{\tau=-\nj{j}+1}^{\nj{j}-1} (\nj{j}-|\tau|)
\PCov^2(\overline{\dwt}_{j,0},\overline{\dwt}_{j,\tau
}) <
\infty.
\end{equation}
Now, consider $\{ \tilde{\dwt}_{j,k} \}$. The Cauchy--Schwarz inequality
and the stationarity of the within-scale process imply that
$\PCov^2(\tilde{\dwt}_{j,0},\tilde{\dwt}_{j,\tau})\leq\PVar
^2(\tilde
{\dwt}_{j,0})=\sigma^{4}_{j}(\tilde{\nu})=O(2^{2j(1-2\alpha)})$,
by~(\ref{eq:vjtildefBound}), and since $\sigma^{2}_{j}(\nu)\asymp
2^{2j\td
}$, we get
%
\begin{equation}
\label{eq:cuvSqrBound2}
\hspace*{6pt}\sup_{n \geq1} \sup_{j=1,\dots,\nmaxscale}
\frac{2^{2j(2\alpha+2\td-1)}}{\nj{j}^2\sigma^{4}_{j}(\nu)}\sum_{\tau
=-\nj
{j}+1}^{\nj{j}-1} (\nj{j} - |\tau|) \PCov^2(\tilde{\dwt
}_{j,0},\tilde
{\dwt}_{j,\tau}) < \infty.
\end{equation}
Finally, using the fact that, for any $j\geq1$, $\bdens{j,0}{\lambda}{f}$ is the
spectral density of the process
$\{\overline{\dwt}_{j,k}\}$ and denoting by $\tilde{\mes}_j$ the
spectral measure of $\{\tilde{\dwt}_{j,k}\}_{k\in\Zset}$, it is
straightforward to show that
\begin{eqnarray*}
A(n,j)&\eqdef&\sum_{\tau=-\nj{j}+1}^{\nj{j}-1}(\nj{j}-|\tau
|)
\PCov(\overline{\dwt}_{j,0},\overline{\dwt}_{j,\tau})\PCov
(\tilde{\dwt
}_{j,0},\tilde{\dwt}_{j,\tau})\\
&=&\int_{-\pi}^\pi\int_{-\pi}^\pi\bdens{j,0}{\lambda'}{f}
\Biggl|\sum_{k=0}^{\nj{j}-1}\rme^{\rmi k(\lambda+\lambda')}\Biggr|^2
\tilde{\mes}_j(d\lambda)\,d\lambda'\\
&\leq&2\pi\nj{j} \sigma^{2}_{j}(\tilde{\nu}) \|\bdens{j,0}{\cdot}{f}\|_\infty .
\end{eqnarray*}
This implies that $A(n,j)\geq0$ and using~(\ref{eq:vjtildefBound}),~(\ref
{eq:DjApproxbis}) and
$\sigma^{2}_{j}(\nu)\asymp2^{2j\td}$, we get
%
\begin{equation}
\label{eq:doubleproduiWbarWtilde}
\sup_{n \geq1} \sup_{j=1,\dots,\nmaxscale}\frac{2^{j(2\alpha
+2\td
-1)}}{\nj{j}\sigma^{4}_{j}(\nu)}
|A(n,j)
| <\infty.
\end{equation}
Using the fact that $\dwt_{j,k}=\overline{\dwt}_{j,k}+\tilde{\dwt}_{j,k}$ and
$\overline{\dwt}_{j,k}$ and $\tilde{\dwt}_{j,k}$
are uncorrelated,~(\ref{eq:emvarLinCase}), (\ref{eq:cumfourBound}),
(\ref{eq:cuvSqrBound1}), (\ref{eq:cuvSqrBound2}),
~(\ref{eq:doubleproduiWbarWtilde}) and $1-2\alpha-2\td<-\beta<0$
yield~(\ref{eq:rosenthalAssump}).
\end{pf*}
\begin{remark}\label{remvarepsilonispi}
If $\varepsilon=\pi$ in the assumptions of Theorem~\ref
{thm:AssumpCons}, then, in the above proof, $\tilde{\dwt}_{j,k}=0$ for
all $(j,k)$, so not only~(\ref{eq:rosenthalAssump}) holds, but also
the stronger relation
%
\begin{equation}
\label{eq:rosenthalAssumpBis}
\sup_{n \geq1} \sup_{j=1,\dots,\nmaxscale} \nj{j}^{-1} \PVar
\Biggl(\sum
_{k=0}^{\nj{j}-1} \frac{\dwt_{j,k}^{2}}{\sigma^{2}_{j}(\nu)}\Biggr) <
\infty .
\end{equation}
\end{remark}
\begin{pf*}{Proof of Corollary~\protect\ref{cor:C1gaussian}}
Condition~\ref{cond:GeneralCondition} holds because
Theorem~\ref{thm:AssumpCons} applies to a Gaussian process.
Moreover, since its fourth order
cumulants are zero, the relation
$\dwt_{j,k}^2=\overline{\dwt}_{j,k}^2+\tilde{\dwt
}_{j,k}^2+2\overline
{\dwt}_{j,k}\tilde{\dwt}_{j,k}$,~(\ref{eq:cuvSqrBound2})
and~(\ref{eq:doubleproduiWbarWtilde}) yield
\[
\PVar\Biggl(\sum_{k=0}^{\nj{j}-1} (\dwt_{j,k}^{2}-\overline{\dwt
}_{j,k}^{2})\Biggr)
\leq C
\biggl[\frac{\nj{j}^2\sigma^{4}_{j}(\nu)}{2^{2j(2\alpha+2\td
-1)}}+\frac{\nj
{j}\sigma^{4}_{j}(\nu)}{2^{j(2\alpha+2\td-1)}}\biggr],
\]
where $C$ is a positive constant. Since $\overline{\dwt}_{j,k}$
and $\tilde{\dwt}_{j,k}$ are uncorrelated, $\PE[\dwt
_{j,k}^{2}-\overline
{\dwt}_{j,k}^{2}]=\sigma^{2}_{j}(\tilde{\nu})$, hence the
last display, $\sigma^{2}_{j}(\nu)\asymp2^{2j\td}$ and~(\ref
{eq:vjtildefBound}) yield~(\ref{eq:gaussianCaseApproxXbar}).
\end{pf*}

\section{Asymptotic behavior of the contrast process}\label
{sec:asympt-behav-contr}
We decompose the contrast~(\ref{eq:TildeJdef}) into a sum of a
(deterministic) function of $d$ and a random process indexed by $d$,
%
\begin{equation}
\label{eq:Contrast2}
\contrast_{\indexset} (d) \eqdef\contrastdet_{\indexset}(d) +
\fluctuation_{\indexset}(d) +
\log\bigl(\cardinal{\indexset} \K2^{2\td\jmean{\indexset
}}\bigr),
\end{equation}
where the log term does not depend on $d$ (and thus may be discarded) and
\begin{eqnarray}
\contrastdet_{\indexset}(d) &\eqdef&\log\Biggl( \frac1{\cardinal
{\indexset}} \sum_{(j,k) \in\indexset} 2^{2(\td-d)j} \Biggr) -
\frac{1}{\cardinal{\indexset}} \sum_{(j,k) \in\indexset} \log
\bigl(2^{2(\td-d)j}\bigr) ,\label{eq:contrast_deterministe} \\
\fluctuation_{\indexset}(d) &\eqdef&\log\Biggl[1 +
\sum_{(j,k) \in\indexset} \frac{2^{2(\td-d)j}}{\sum_{\indexset}
2^{2(\td-d)j}}
\biggl(\frac{\dwt_{j,k}^2}{\K 2^{2 \td j}}-1\biggr)
\Biggr],\label{eq:fluctuation}
\end{eqnarray}
with $\K$ defined in~(\ref{eq:BiasTermAssump}).
\begin{proposition}
\label{prop:DetermPart}
For any finite and nonempty set $\indexset\subset\Nset\times\Zset$,
the function $d \to\contrastdet_{\indexset}(d)$ is nonnegative,
convex and vanishes at $d=\td$.
Moreover, for any sequence $\{\nlowscale\}$ such that $n2^{-\nlowscale}
\to\infty$ as $n\to\infty$, and for any constants
$\admin$ and $\admax$ in $\Rset$ satisfying $\td-1/2 < \admin\leq
\admax$,
%
\begin{equation}
\label{eq:AsympEquivalent}
\liminf_{n\to\infty}
\inf_{d \in[\admin,\admax]} \inf_{j_1=\nlowscale+1,\dots
,\nmaxscale
}\ddot{\contrastdet}_{\indexset_n(\nlowscale,j_1)}(d) > 0 ,
\end{equation}
where $\indexset_n$ is defined in~\textup{(\ref{eq:InDef})} and $\ddot
{\contrastdet}_\indexset$ denotes the second derivative of
$\contrastdet_\indexset$.
\end{proposition}
\begin{pf*}{Proof}
By concavity of the $\log$ function, $\contrastdet_\indexset(d) \geq0$
and is zero if $d= \td$.
If $\indexset=\indexset_n(\nlowscale,j_1)$ with $j_1\geq\nlowscale+1$,
one can compute $\ddot{\contrastdet}_\indexset(d)$
and show that it can be expressed as $\ddot{\contrastdet}_\indexset(d)
= (2 \log(2))^2
\PVar(N)$, where $N$ is an integer-valued random variable such that
$ \prob(N=j)=2^{2(\td-d)j}  \nj{j} / \sum_{j= \nlowscale}^{j_1}
2^{2(\td-d)j} \nj{j}$ for $j\geq0$.
Let $d\geq\admin> \td-1/2$. Then,
\[
\prob(N=\nlowscale) 
\geq\bigl(1-2^{2(\td-\admin)-1}\bigr) \{ 1 - {\L}
2^{\nlowscale}(n-\L+1)^{-1}\}.
\]
Since $n2^{-\nlowscale} \to\infty$, the term between the brackets
tends to 1 as $n\to\infty$. Hence, for $n$ large enough, we
have $\inf_{d\geq\admin}\prob(N=\nlowscale)\geq(1-2^{2(\td
-\admin
)-1}) /2$. Similarly, one finds, for $n$ large enough,
$\inf_{d\in[\admin,\admax]}\prob(N=\nlowscale+1)\geq(1-2^{2(\td
-\admin
)-1}) 2^{2(\td-\admax)-1} /2$.
Hence,
\begin{eqnarray*}
\inf_{d\in[\admin,\admax]} \PVar(N) &\geq&\{\nlowscale-\PE(N)\}
^2\prob(N=\nlowscale)\\
&&{}+\{\nlowscale+1-\PE(N)\}^2\prob(N=\nlowscale+1) \\
& \geq&\bigl(1-2^{2(\td-\admin)-1}\bigr) 2^{2(\td-\admax)-2}\\
&&{}\times\bigl(\{
\nlowscale-\PE(N)\}^2+\{\nlowscale+1-\PE(N)\}^2\bigr)\\
&\geq&\bigl(1-2^{2(\td-\admin)-1}\bigr) 2^{2(\td-\admax)-4} ,
\end{eqnarray*}
where the last inequality is obtained by observing that either $\PE
(N)-\nlowscale\geq1/2$ or $\nlowscale+1-\PE(N)<1/2$.
\end{pf*}

We now show that the random component $\fluctuation_{\indexset}(d)$ of
the contrast~(\ref{eq:Contrast2}) tends to~0 uniformly
in $d$.
For all $\rho>0$, $q\geq0$ and $\delta\in\Rset$, define the set of
real-valued sequences
%
\begin{equation}\label{eq:calBdef}
\calB(\rho,q,\delta)\eqdef\bigl\{ \{\mu_j\}_{j\geq0}\dvtx|\mu
_j|\leq
\rho (1+j^q)2^{j\delta}\mbox{ for all $j\geq0$}\bigr\} .
\end{equation}
Define, for any $n\geq1$, any sequence $\bmu\eqdef\{\mu_j\}_{j\geq0}$
and $0\leq j_0\leq j_1\leq\nmaxscale$,
%
\begin{equation}\label{eq:empSbisJ1}
\empSbis_{n,j_0,j_1}(\bmu)\eqdef\sum_{j=j_0}^{j_1} \mu_{j-j_0}
\sum_{k=0}^{\nj{j}-1} \biggl[\frac{\dwt_{j,k}^2}{\K 2^{2 \td
j}}-1\biggr].
\end{equation}
%
%
\begin{proposition}\label{prop:consistencyIngredients2Improved}
Under Condition~\ref{cond:GeneralCondition}, for any
$q\geq0$ and $\delta<1 $, there exists $C>0$ such that for all $\rho
\geq0$, $n\geq1$ and $j_0=1,\dots,\nmaxscale$,
%
\begin{eqnarray}
\label{eq:uniformEmpBound}
&&\hspace*{6pt}\biggl\{\PE\sup_{\bmu\in\calB(\rho,q,\delta)}\sup
_{j_1=j_0,\dots
,\nmaxscale}|\empSbis_{n,j_0,j_1}(\bmu)|^2\biggr\}^{1/2}\nonumber\\
&&\hspace*{6pt}\qquad \leq C \rho n2^{-j_0}  [H_{q,\delta}(n2^{-j_0})+2^{-\beta
j_0}] ,\\
&&\hspace*{6pt}\mbox{where, for all $x\geq0$, }
H_{q,\delta}(x)\eqdef
\cases{
x^{-1/2}, &\quad$\mbox{if }\delta<1/2 $,\cr
\log^{q+1}(2+x) x^{-1/2}, &\quad$\mbox{if }\delta=1/2 $,\cr
\log^{q}(2+x)x^{\delta-1},&\quad$\mbox{if }\delta>1/2$.
}\nonumber
\end{eqnarray}
\end{proposition}
\begin{pf*}{Proof}
We set $\rho=1$ without loss of generality. We write
\[
\empSbis_{n,j_0,j_1}(\bmu) =
\sum_{j=j_0}^{j_1}\frac{\sigma^{2}_{j}(\nu)}{\K\,2^{2 \td j}} \mu
_{j-j_0}
\sum_{k=0}^{\nj{j}-1}\biggl[\frac{\dwt_{j,k}^2}{\sigma^{2}_{j}(\nu)}-1\biggr]
+ \sum_{j=j_0}^{j_1} \nj{j} \mu_{j-j_0}\biggl[\frac{\sigma^{2}_{j}(\nu)}{\K
 2^{2 \td j}}-1\biggr]
\]
and denote the two terms of the right-hand side of this equality as
$\empSbis^{(0)}_{n,j_0,j_1}(\bmu)$ and $\empSbis
^{(1)}_{n,j_0,j_1}(\bmu
)$, respectively.
By~(\ref{eq:BiasTermAssump}), $C_1\eqdef\sup_{j \geq0} 2^{\beta j}
|\sigma^{2}_{j}(\nu)/(\K 2^{2 \td j})-1| < \infty$,
which implies $\sup_{j\geq0}|\sigma^{2}_{j}(\nu)/(\K 2^{2 \td
j})|\leq
1+C_1$. Hence, if $\bmu\in\calB(1,q,\delta)$, then
\[
\bigl|\empSbis^{(0)}_{n,j_0,j_1}(\bmu)\bigr|
\leq(1+C_1) \sum_{j=j_0}^{\nmaxscale}\bigl(1+(j-j_0)^q\bigr)
2^{(j-j_0)\delta}
\Biggl|\sum_{k=0}^{\nj{j}-1}\biggl(\frac{\dwt_{j,k}^2}{\sigma^{2}_{j}(\nu)}-1\biggr)\Biggr|.
\]
Using the Minkowski inequality and $\nj{j}\leq n2^{-j}$, (\ref
{eq:rosenthalAssump}) implies that there exists a constant
$C_2$ such that
\begin{eqnarray}\label{eq:Dtildedtilde}
\hspace*{12pt}&&\biggl\{\PE\biggl[ \sup_{\bmu\in\calB(1,q,\delta)}\sup
_{j_1=j_0,\dots
,\nmaxscale}
\bigl|\empSbis^{(0)}_{n,j_0,j_1}(\bmu)\bigr|^2\biggr]\biggr\}
^{1/2} \nonumber\\[-8pt]\\[-8pt]
&&\qquad \leq(1+C_1) C_2 \sum_{j=j_0}^{\nmaxscale}\bigl(1+(j-j_0)^q\bigr)
2^{(j-j_0)\delta}\bigl[(n2^{-j})^{1/2}+n2^{-(1+\beta)j}\bigr].\nonumber
\end{eqnarray}
The sum over the first term is $O(n2^{-j_0}H_{q,\delta
}(n2^{-j_0}))$ since
$\nmaxscale-j_0\asymp\log_2n+\log_22^{-j_0}=\log_2(n2^{-j_0})$.
The sum
over the second term is $O(n2^{-(1+\beta)j_0})$ since
$\delta<1$ and $1+\beta>1$, so (\ref{eq:Dtildedtilde}) is
$O((n2^{-j_0})\{H_{q,\delta}(n2^{-j_0})+2^{-\beta
j_0}\})$ since $2^{\nmaxscale}\asymp n$. Now, by the definition of
$C_1$ above and since $\nj{j}\leq n2^{-j}$, we get
\[
\sup_{\bmu\in\calB(1,q,\delta)}\sup_{j_1=j_0,\dots,\nmaxscale
}\bigl|
\empSbis^{(1)}_{n,j_0,j_1}(\bmu)\bigr|
\leq C_1 n   \sum_{j=j_0}^{\nmaxscale}\bigl(1+(j-j_0)^q\bigr)
2^{(j-j_0)\delta
}2^{-j(1+ \beta)} ,
\]
which is $O(n2^{-(1+\beta) j_0})$. The two last displays
yield~(\ref{eq:uniformEmpBound}).
\end{pf*}
\begin{corollary}\label{cor:calEb}
Let $\{\nlowscale\}$ be a sequence such that $\nlowscale
^{-1}+(n2^{-\nlowscale})^{-1}\to0$ as $n\to\infty$ and let
$\fluctuation_{\indexset}(d)$ be defined as in~\textup{(\ref{eq:fluctuation})}.
Condition~\ref{cond:GeneralCondition} then implies that as
$n\to\infty$:
\begin{longlist}[(a)]
\item[(a)]\label{itemCorCalEb-a} for any $\k\geq0$,
\[
\sup_{d\in\Rset} \bigl|\fluctuation_{\indexset_n(\nlowscale
,\nlowscale
+\k)}(d)\bigr|
=O_{\prob}\bigl( (n2^{-\nlowscale})^{-1/2}+2^{-\beta\nlowscale
}\bigr);
\]
\item[(b)]\label{itemCorCalEb-b} for all $\admin>\td-1/2$, setting
$\delta
=2(\td-\admin)$,
\[
\sup_{d\geq\admin}\sup_{j_1=\nlowscale,\dots,\nmaxscale}
\bigl|\fluctuation_{\indexset_n(\nlowscale,j_1)}(d)\bigr|
=O_{\prob}\bigl( H_{0,\delta}(n2^{-\nlowscale})+2^{-\beta
\nlowscale
}\bigr) .
\]
\end{longlist}
\end{corollary}
\begin{pf*}{Proof}
The definitions \eqref{eq:fluctuation} and~ \eqref{eq:empSbisJ1} imply
that, for $0\leq j_0 \leq j_1\leq\nmaxscale$,
\[
\fluctuation_{\indexset_n(j_0,j_1)}(d)=\log\bigl[1+
(n2^{-j_0})^{-1}\empSbis_{n,j_0,j_1}[\bmu(d,j_0,j_1)]\bigr]
\]
with $\bmu(d,j_0,j_1)$ is the sequence $\{\mu_{j}(d,j_0,j_1)\}_{j\geq
0}$ defined by
\begin{equation}\label{eq:mujdDef}
\mu_{j}(d,j_0,j_1)  \eqdef n2^{-j_0}   \frac{2^{2 (\td-d) (j+j_0)}}
{\sum_{j'=j_0}^{j_1} 2^{2 (\td-d)j'}\nj{j'}}   \1(j \leq j_1-j_0).
\end{equation}
The bounds \hyperref[itemCorCalEb-a]{(a)} and \hyperref[itemCorCalEb-b]{(b)} then follow
from Proposition~\ref{prop:consistencyIngredients2Improved}, the
Markov inequality and the following bounds.

\textit{Part} \hyperref[itemCorCalEb-a]{(a)}.
In this case, we apply Proposition~\ref{prop:consistencyIngredients2Improved} with $\delta=0$.
Indeed, using the fact that $\mu_j(d,\nlowscale,\nlowscale+\k)\nj
{\nlowscale
+j}\leq n2^{-\nlowscale}$ for all $j=0,\dots,\k$ and is zero otherwise,
we have that $\mu_j\leq n2^{-\nlowscale}/\nj{\nlowscale+\k}\to
2^\k$
as $n\to\infty$, since $n2^{-\nlowscale}\to\infty$. Then, for
large enough
$n$, $\bmu(d,\nlowscale,\nlowscale+\k)\in\calB(2^{\k+1},0,0)$
for all
$d\in\Rset$.

\textit{Part} \hyperref[itemCorCalEb-b]{(b)}. Here, we still
apply Proposition~\ref{prop:consistencyIngredients2Improved}, but with
$\delta=2(\td-\admin)<1$, implying that
$H_{0,\delta}(n2^{-\nlowscale})\to0$.
Indeed, since the denominator of the ratio appearing in~(\ref{eq:mujdDef})
is at least $2^{2 (\td-d)\nlowscale}\nj{\nlowscale}$,
we have
$\sup_{j_1\geq\nlowscale}\sup_{d\geq\admin}|\mu_j(d,\break\nlowscale,j_1)|
\leq n2^{-\nlowscale} \nj{\nlowscale}^{-1} 2^{\delta
j}$. Since $n2^{-\nlowscale}\sim\nj{\nlowscale}$ as $n\to\infty$, we
get that,
for large enough $n$, $\bmu(d,\nlowscale,j_1) \in\calB(2,0,\delta
)$ for
all $d\geq\admin$ and $j_1\geq\nlowscale$.
\end{pf*}

\section{Weak consistency}\label{sec:weak-consistency}
We now establish a preliminary result on the consistency of $\hd$. It
does not provide an optimal rate, but it will be used in
the proof of Theorem~\ref{thm:Rates}, which provides the optimal rate.
By the definition of $\hd$ and \eqref{eq:Contrast2}, we have
%
\begin{equation}\label{eq:ContrastTriangIneq}
0\geq\contrast_\indexset(\hd_\indexset) - \contrast_\indexset
(\td)
= \contrastdet_\indexset(\hd_\indexset)+ \fluctuation_\indexset
(\hd
_\indexset) - \fluctuation_\indexset(\td).
\end{equation}
The basic idea for proving consistency is to show that (1) the function
$d \mapsto\contrast(d)$ behaves as $(d- \td)^2$ up to a multiplicative
positive constant and (2) the function $d \mapsto\fluctuation(d)$ tends
to zero in probability, uniformly in $d$.
Proposition~\ref{prop:DetermPart} will prove (1) and Corollary~\ref{cor:calEb} will yield (2).
\begin{proposition}[(Weak consistency)]
\label{prop:consistencyGaussian}
Let $\{\nlowscale\}$ be a sequence such that $\nlowscale
^{-1}+(n2^{-\nlowscale})^{-1}\to0$ as $n\to\infty$.
Condition~\textup{\ref{cond:GeneralCondition}} implies that as $n\to\infty$,
%
\begin{equation}\label{eq:WCvitessefixed}
\sup_{j_1=\nlowscale+1,\dots,\nmaxscale}\bigl|\hd_{\indexset
_n(\nlowscale,j_1)}-\td\bigr| = O_{\prob}\{
(n2^{-\nlowscale})^{-1/4}+2^{-\beta\nlowscale/2} \} .
\end{equation}
\end{proposition}
\begin{pf*}{Proof}
The proof proceeds in four steps.
\begin{enumerate}[\textit{Step} 1.]
\item[\textit{Step} 1.]\label{step:0} For any positive integer $\k$,
$|\hd_{\indexset_n(\nlowscale,\nlowscale+\k)}-\td| =
o_{\prob
}(1) .$
\item[\textit{Step} 2.]\label{step:1} There exists $d_{\min}\in(\td-1/2,\td)$ such that,
as $n\to\infty$,
\[
\prob\biggl\{\inf_{j_1=\nlowscale+2,\dots,\nmaxscale}\hd
_{\indexset
_n(\nlowscale,j_1)}\leq d_{\min}\biggr\} \to0 .
\]
Combining this with~Step~\textup{\hyperref[step:0]{1}} yields
$\prob\{\inf_{j_1=\nlowscale+1,\dots,\nmaxscale
}\hd_{\indexset_n(\nlowscale,j_1)}\leq d_{\min}\}\break \to~0$.
\item[\textit{Step} 3.]\label{step:2} For any $d_{\max}>\td$, as $n\to\infty$,
$\prob\{\sup_{j_1=\nlowscale+1,\dots,\nmaxscale
}\hd
_{\indexset_n(\nlowscale,j_1)}\geq d_{\max}\} \to0$.
\item[\textit{Step} 4.]\label{step:3} Define $H_{0,\delta}$ as in
Proposition~\ref{prop:consistencyIngredients2Improved}. For all $\admin\in(\td
-1/2,\td
)$ and $\admax>\td$, setting $\delta=2(\td-\admin)$, we have
\begin{eqnarray*}
&&\sup_{j_1=\nlowscale+1,\dots,\nmaxscale}\bigl[\1_{[\admin,\admax
]}\bigl(\hd_{\indexset_n(\nlowscale,j_1)}\bigr) \bigl(\hd_{\indexset_n(\nlowscale
,j_1)}-\td
\bigr)^2\bigr]\\
&&\qquad = O_{\prob}\bigl( H_{0,\delta}(n2^{-\nlowscale}) + 2^{-\beta
\nlowscale
} \bigr).
\end{eqnarray*}
\end{enumerate}
Before proving these four steps, let us briefly explain how they
yield~(\ref{eq:WCvitessefixed}).
First, observe that they imply that
$\sup_{j_1=\nlowscale+1,\dots,\nmaxscale}|\hd_{\indexset
_n(\nlowscale
,j_1)}-\td|=o_\prob(1)$.
Then, applying~Step~\textup{\hyperref[step:3]{4}} again with $\admin\in(\td-1/4,\td)$, so
that $H_{0,\delta}(x)=x^{-1/2}$, we obtain~(\ref{eq:WCvitessefixed}).
\end{pf*}
\begin{pf*}{Proof of Step \textup{\protect\hyperref[step:0]{1}}}
Using standard arguments for contrast estimation (similar to those
detailed is Step~\textup{\hyperref[step:2]{3}} and Step~\textup{\hyperref[step:3]{4}} below), this
step is a direct consequence of Proposition~\ref{prop:DetermPart} and
Corollary~\ref{cor:calEb}(\hyperref[itemCorCalEb-a]{a}).
\end{pf*}
\begin{pf*}{Proof of Step~\textup{\protect\hyperref[step:1]{2}}}
Using~(\ref{eq:TildeJdef}), we have, for all $d\in\Rset$,
\[
\contrast_\indexset(d) - \contrast_\indexset(\td)
=\log\Biggl(\sum_{(j,k)\in\indexset} 2^{2(d-\td)(\jmean{\indexset}-j)}
\frac{\dwt_{j,k}^2}{\K  2^{2 \td j}}\Biggr)
- \log\Biggl(\sum_{(j,k)\in\indexset} \frac{\dwt_{j,k}^2}{\K
2^{2 \td
j}} \Biggr).
\]
For some $d_{\min} \in(\td-1/2,\td)$ to be specified later, we set
%
\begin{equation}\label{eq:wjd}
\quad w_{\indexset,j}(d) \eqdef2^{2(j-\jmean{\indexset})(\td-d)} \1 \{
j\leq
\jmean{\indexset} \} + 2^{2(j-\jmean{\indexset}) (\td- d_{\min})}
\1 \{
j > \jmean{\indexset} \} ,
\end{equation}
so that for all $j$ and $d\leq d_{\min}$, $w_{\indexset,j}(d)\leq
2^{2(d-\td)(\jmean{\indexset}-j)}$. We further obtain, for all
$d\leq
d_{\min}$,
%
\begin{equation}\label{eq:ContrastIneq1}
\contrast_\indexset(d) - \contrast_\indexset(\td)
\geq\log\frac{\Sigma_\indexset(d) + A_\indexset(d)}{1 +
B_\indexset} ,
\end{equation}
where $\Sigma_\indexset(d)\eqdef\cardinal{\indexset}^{-1} \sum
_{(j,k)\in\indexset} w_{\indexset,j}(d)$, $A_\indexset(d) \eqdef
\cardinal{\indexset}^{-1} \sum_{(j,k)\in\indexset}
w_{\indexset,j}(d)\times\break
(\frac{\dwt_{j,k}^2}{\K  2^{2 \td j}}-1)$ and
$B_\indexset\eqdef\cardinal{\indexset}^{-1} \sum_{(j,k)\in
\indexset}
(\frac{\dwt_{j,k}^2}{\K  2^{2 \td j}}-1)$.
We will show that $d_{\min} \in(\td-1/2,\td)$ may be chosen in such a
way that
\begin{eqnarray}
\label{eq:SigmaBound}
\liminf_{n \to\infty} \inf_{d \leq d_{\min}} \inf_{j_1=\nlowscale
+2,\dots,\nmaxscale} \Sigma_{\indexset_n(\nlowscale,j_1)}(d) &>& 1, \\
\label{eq:AandBBound}
\sup_{j_1=\nlowscale+2,\dots,\nmaxscale} \biggl(\sup_{d\leq
d_{\min}}
\bigl|A_{\indexset_n(\nlowscale,j_1)}(d)\bigr|
+\bigl|B_{\indexset_n(\nlowscale,j_1)}\bigr|\biggr)&=& o_{\prob}(1)  .
\end{eqnarray}
By (\ref{eq:ContrastTriangIneq}), $\contrast_\indexset(\hd
_\indexset)
\leq\contrast_\indexset(\td)$. Then,
$\inf_{j_1=\nlowscale+2,\dots,\nmaxscale} \hd_{\indexset
_n(\nlowscale
,j_1)} \leq d_{\min}$
would imply that there exists $j_1=\nlowscale+2,\dots,\nmaxscale$ such
that $\inf_{d \leq d_{\min}}
\contrast_{\indexset_n(\nlowscale,j_1)}(d)-\contrast(\td)\leq0$, an
event whose probability tends to zero as a consequence
of~(\ref{eq:ContrastIneq1})--(\ref{eq:AandBBound}). Hence, these
equations yield~Step~\textup{\hyperref[step:1]{2}}.
It thus remains to show that~(\ref{eq:SigmaBound})
and~(\ref{eq:AandBBound}) hold. By Lemma~\ref{lem:etavareta},
since $n2^{-\nlowscale}\to\infty$, we have, for $n$ large
enough,
%
\begin{equation}\label{eq:etaJ0}
\sup_{j_1=\nlowscale,\dots,\nmaxscale}\jmean{\indexset
_n(\nlowscale
,j_1)}< \nlowscale+1  .
\end{equation}
Using $w_{\indexset_n(\nlowscale,j_1),\nlowscale}(d)\geq0$ and, for $n$
large enough, $w_{\indexset_n(\nlowscale,j_1),j}(d)\geq\break
2^{2(j-(\nlowscale+1)) (\td- d_{\min})} $,
for $j\geq\nlowscale+1$, we get, for all $d\leq d_{\min}<\td$ and
$j_1=\nlowscale+2,\ldots,\nmaxscale$,
\[
\Sigma_{\indexset_n(\nlowscale,j_1)}(d)
\geq\frac{2^{-2(\nlowscale+1) (\td-d_{\min})}}{\cardinal
{\indexset
_n(\nlowscale,\nmaxscale)}} \sum_{j=\nlowscale+1}^{\nlowscale+2}
2^{2j(\td-d_{\min})} \nj{j}  .
\]
Since $n2^{-\nlowscale}\to\infty$, using Lemma~\ref{lem:etavareta},
$n\asymp2^{\nmaxscale}$ and the fact that $2(\td-d_{\min})-1<0$,
straightforward computations give that the LHS in the
previous display is asymptotically equivalent to
$(1-2^{\{2(\td-d_{\min})-1\}2})/(4-2^{2(\td-d_{\min})+1})$.
There are values of $d_{\min} \in(\td-1/2,\td)$
such that this ratio is strictly larger than 1. For such a choice and
for $n$ large enough,~(\ref{eq:SigmaBound}) holds.

We now check~(\ref{eq:AandBBound}). Observing that, for $\indexset_n
\eqdef\indexset_n(\nlowscale,j_1)$ and using the
notation~(\ref{eq:empSbisJ1}),
$A_{\indexset_n}(d)= \cardinal{\indexset_n}^{-1} |\empSbis
_{n,\nlowscale
,j_1}( \{w_{\indexset_n,\nlowscale+j}(d)\})|$ and $B_{\indexset_n} =
\cardinal{\indexset_n}^{-1}
\empSbis_{n,\nlowscale,j_1}(\1)$, the bound~(\ref{eq:AandBBound})
follows from $\cardinal{\indexset_n}\geq\nj{\nlowscale}\sim
n2^{-\nlowscale}$ and
Proposition~\ref{prop:consistencyIngredients2Improved} since, for all
$d \leq d_{\min}$ and $j\geq0$,
$w_{\indexset_n,\nlowscale+j}(d) \leq2^{2(\nlowscale+j-\jmean
{\indexset
_n})(\td-d_{\min})}\leq2^{2j(\td-d_{\min})}$, which shows that
$\{w_{\indexset_n,\nlowscale+j}(d)\}_{j\geq0}$ belongs to $\calB
(1,0,\delta)$ with $\delta=2(\td-d_{\min}))<1$.
\end{pf*}
\begin{pf*}{Proof of Step \textup{\protect\hyperref[step:2]{3}}}
By~(\ref{eq:ContrastTriangIneq}), $\contrastdet_\indexset(\hd
_\indexset
) \leq\fluctuation_\indexset(\td) -\fluctuation_\indexset(\hd
_\indexset)$,
so, for any $d_{\max} \geq\td$, one has $\inf_{d\geq d_{\max}}
\contrastdet_\indexset(d) \leq2\sup_{d\geq\td}|
\fluctuation_\indexset(d) |$ on the event $\{ \hd_\indexset
\geq
d_{\max} \}$.
By Proposition~\ref{prop:DetermPart}, there exists $c>0$ such that, for
$n$ large enough,
$\contrastdet_{\indexset_n(\nlowscale,j_1)}(d)\geq c$ uniformly for $d
\geq d_{\max}$ and
$j_1=\nlowscale+1,\dots,\nmaxscale$. Thus, for $n$ large enough,
\[
\prob\biggl\{ \sup_{j_1=\nlowscale+1,\dots,\nmaxscale}\hd
_{\indexset
_n(\nlowscale,j_1)} \geq d_{\max} \biggr\}
\leq\prob\biggl\{ 2 \sup_{ d \geq\td}\sup_{j_1=\nlowscale
+1,\dots,\nmaxscale} \bigl|\fluctuation_{\indexset_n(\nlowscale,j_1)}(d)\bigr| \geq
c\biggr\},
\]
which tends to 0 as $n\to\infty$, by Corollary~\ref{cor:calEb}(\hyperref[itemCorCalEb-b]{b}).
\end{pf*}
\begin{pf*}{Proof of Step \textup{\protect\hyperref[step:3]{4}}}
Equation \eqref{eq:ContrastTriangIneq} implies that
$\1_{[\admin,\admax]}(\hd_\indexset) \contrastdet_\indexset(\hd
_\indexset)\leq2 \sup_{d\geq\admin} | \fluctuation
_\indexset(d)|$.
Let $c$ denote the liminf in the left-hand side
of~(\ref{eq:AsympEquivalent}) when $\admin=\admin$ and $\admax=\admax$.
Proposition \ref{prop:DetermPart} and a second order Taylor expansion
of $\contrastdet_\indexset$ around $\td$ give that,
for $n$ large enough, for all $j_1=\nlowscale+1,\dots,\nmaxscale$ and
$d\in[\admin,\admax]$,
$\contrastdet_{\indexset_n(\nlowscale,j_1)}(d)\geq(c/4) (d-\td)^2$.
Hence, for $n$ large enough,\looseness=1
\[
\sup_{j_1=\nlowscale+1,\dots,\nmaxscale}\bigl[\1_{[\admin,\admax
]}\bigl(\hd
_{\indexset_n(\nlowscale,j_1)}\bigr)\bigl(\hd_{\indexset_n(\nlowscale
,j_1)}-\td
\bigr)^2 \bigr]
\leq\frac8c  \sup_{d\geq\admin} \bigl| \fluctuation_{\indexset
_n(\nlowscale,j_1)}(d) \bigr| .\vadjust{\goodbreak}
\]
Corollary~\ref{cor:calEb}(\hyperref[itemCorCalEb-b]{b}) then yields Step~\textup{\hyperref[step:3]{4}}.
\end{pf*}
\begin{remark}
Proposition~\ref{prop:consistencyGaussian} implies that if
$\nlowscale
\leq\nupscale\leq\nmaxscale$ with
$\nlowscale^{-1}+(n2^{-\nlowscale})^{-1}\to0$ as $n\to\infty$,
then $\hd
_{\indexset_n(\nlowscale,\nupscale)}$ is a consistent estimator of
$\td$.
While the rate provided by~(\ref{eq:WCvitessefixed}) is not optimal, it
will be used to derive
the optimal rates of convergence (Theorem~\ref{thm:Rates}).
\end{remark}

\section{\texorpdfstring{Proofs of Theorems~\protect\ref{thm:Rates} and \protect\ref{thm:CLT}}{Proofs of Theorems~3 and 5}}\label{sec:rates}

\begin{nation*}
In the following, $\{\nlowscale\}$ and $\{\nupscale\}$ are two
sequences satisfying~(\ref{eq:J0J1n}).
The only difference between the two following settings
(\hyperref[it:caseJ1-J0fixed]{S-1})
(where $\nupscale-\nlowscale$ is fixed) and (\hyperref[it:caseJ1isJ]{S-2}) (where
$\nupscale-\nlowscale\to\infty$) lies in the computations of the
asymptotic variances in
Theorem~\ref{thm:CLT} (CLT). Hence, we shall hereafter write
$\lowscale$, $\upscale$, $\indexset_n$, $\hd_n$, $\Sclt_n$ and
$\empSbis_n$ 
for $\nlowscale$, $\nupscale$, $\indexset_n(\nlowscale,\nupscale
)$, $\hd
_{\indexset_n(\nlowscale,\nupscale)}$,
$\Sclt_{\indexset_n(\nlowscale,\nupscale)}$ and $\empSbis
_{n,\nlowscale
,\nupscale}$, respectively.

We will use the explicit notation when the distinction between these
two cases~\textup{(\hyperref[it:caseJ1-J0fixed]{S-1})} and~\textup{(\hyperref[it:caseJ1isJ]{S-2})}
is necessary, namely, when computing the limiting variances in the proof
of Theorem~\ref{thm:CLT}.
\end{nation*}
\begin{pf*}{Proof of Theorem \protect\ref{thm:Rates}}
Since $\Sclt_n(\hd_n)=0$ [see~\eqref{eq:ScltDef}] a Taylor expansion of
$\Sclt_n$ around $d=\hd_n$ yields
%
\begin{equation}
\label{eq:EstEqTaylor}
\Sclt_n(\td) =
2\log(2)  (\hd_n-\td)  \sum_{(j,k)\in\indexset_n} (j-\jmean
{\indexset
_n}) j  2^{-2j \htd_n } \dwt_{j,k}^2
\end{equation}
for some $\htd_n$ between $\td$ and $\hd_n$. The proof of
Theorem~\ref{thm:Rates} now consists of
bounding $\Sclt_n(\td)$ from above and showing that
$\sum_{\indexset_n}   (j-\jmean{\indexset_n}) j  2^{-2j \htd_n}
\dwt
_{j,k}^2$, appropriately normalized, has a strictly positive limit.

By the definitions of $\Sclt_n$ [see~(\ref{eq:ScltDef})], $\empSbis_n$
[see~(\ref{eq:empSbisJ1})] and $\jmean{\indexset_n}$
[see~(\ref{eq:etaDelta})], we have
$ \Sclt_n(\td) = \empSbis_n( \K\{j+\lowscale-\jmean
{\indexset_n}\}
_{j\geq0} )$. Since $\lowscale\leq
\jmean{\indexset_n}\leq\lowscale+1$ for $n$
large enough [see~(\ref{eq:etaJ0})] the sequence $\K\{j+\lowscale
-\jmean{\indexset_n}\}_{j\geq0}$ belongs to $\calB(\K,1,0)$
[see~(\ref{eq:calBdef})], and
Proposition~\ref{prop:consistencyIngredients2Improved}, together with the Markov
inequality, yields, as $n\to\infty$,
%
\begin{eqnarray}\label{eq:ScltBound}
\Sclt_n(\td)&=& n 2^{-\lowscale}O_{\prob}\bigl(H_{1,0}(n
2^{-\lowscale
})+2^{-\beta\lowscale}\bigr)\nonumber\\[-8pt]\\[-8pt]
&=&n 2^{-\lowscale}O_{\prob}\bigl( (n2^{-\lowscale})^{-1/2}+
2^{-\beta
\lowscale} \bigr)  ,\nonumber
\end{eqnarray}
which is the desired upper bound.

We shall now show that the sum in~(\ref{eq:EstEqTaylor}) multiplied by
$n 2^{-\lowscale}$ has a strictly positive lower bound.
Applying Proposition~\ref{prop:consistencyGaussian}, we have
\[
|\htd_n- \td|\leq|\hd_n-\td|=O_{\prob}\bigl( (n2^{-\lowscale})^{-1/4}
+ 2^{-\beta\lowscale/2} \bigr) .
\]
Using the fact that $|2^{2j(\td-\htd_n)}-1|\leq2^{2j|\td-\htd_n|}-1\leq
2\log
(2)j|\td-\htd_n|2^{2j|\td-\htd_n|}$, we have that, on
the event $\{|\td-\htd_n|\leq1/4\}$,
\begin{eqnarray*}
&&\Biggl|\sum_{(j,k)\in\indexset_n}(j-\jmean{\indexset_n}) j
\frac{\dwt
_{j,k}^2}{2^{2\htd_n j}}
- \sum_{(j,k)\in\indexset_n}(j-\jmean{\indexset_n}) j  \frac
{\dwt
_{j,k}^2}{2^{2\td j}} \Biggr|\\
&&\qquad\leq
2\log(2) |\td-\htd_n|   2^{2\lowscale|\td-\htd_n|} \sum
_{(j,k)\in
\indexset_n} |j-\jmean{\indexset_n}| j^2 \frac{\dwt
_{j,k}^2}{2^{2\td
j}} 2^{(j-\lowscale)/2}.
\end{eqnarray*}
Using~(\ref{eq:BiasTermAssump}),~(\ref{eq:etaJ0}), $j^2=(j-\lowscale
)^2+2(j-\lowscale)\lowscale+\lowscale^2$ and $\nj{j}\leq
n2^{-j}$, there is a constant $C>0$ such that
\begin{eqnarray*}
&&\PE\sum_{(j,k)\in\indexset_n} |j-\jmean{\indexset_n}| j^2 \frac
{\dwt
_{j,k}^2}{2^{2\td j}} 2^{(j-\lowscale)/2}\\
&&\qquad \leq Cn2^{-\lowscale}\sum_{j=\lowscale}^\upscale|j-\jmean
{\indexset
_n}| j^22^{-(j-\lowscale)/2}
=O(\lowscale^2n 2^{-\lowscale}).
\end{eqnarray*}
Hence, since $\lowscale^2(n2^{-\lowscale})^{-1/4}\to0$, the last three
displays yield, as $n\to\infty$,
%
\begin{equation}
\label{eq:htdVStd}
\Biggl|\sum_{(j,k)\in\indexset_n}(j-\jmean{\indexset_n}) j
\frac{\dwt
_{j,k}^2}{2^{2\htd_n j}}
- \sum_{(j,k)\in\indexset_n}(j-\jmean{\indexset_n}) j  \frac
{\dwt
_{j,k}^2}{2^{2\td j}} \Biggr|
=o_\prob(n 2^{-\lowscale})  .
\end{equation}
We now write
\begin{eqnarray*}
&&\sum_{(j,k)\in\indexset_n}(j-\jmean{\indexset_n}) j  \frac{\dwt
_{j,k}^2}{2^{2\td j}}\\
&&\qquad =\K\sum_{(j,k)\in\indexset_n}(j-\jmean{\indexset_n}) j
\biggl(\frac
{\dwt_{j,k}^2}{\K2^{2\td j}}-1\biggr)
+ \K\sum_{(j,k)\in\indexset_n}(j-\jmean{\indexset_n}) j .
\end{eqnarray*}
With the notation~(\ref{eq:empSbisJ1}), the first term on the
right-hand side is
$\empSbis_n(\bmu)$, where $\bmu$ is the sequence $\K\{(j+\lowscale
-\jmean{\indexset_n})(j+\lowscale)\}_{j\geq0}$.
In view of~(\ref{eq:etaJ0}), $(j+\lowscale-\jmean{\indexset
_n})(j+\lowscale)\leq j^2+j\lowscale$, so the sequence $\bmu$
is the sum of two sequences belonging to
$\calB(\K,2,0)$ and $\calB(\K\lowscale,1,0)$, respectively. Applying
Proposition~\ref{prop:consistencyIngredients2Improved} together
with the Markov inequality, we get that our $\empSbis_n(\bmu
)=n2^{-\lowscale}O_\prob(H_{0,0}(n2^{-\lowscale})+ \lowscale
H_{1,0}(n2^{-\lowscale}))=n2^{-\lowscale}o_\prob(1)$ since
$\lowscale
(n2^{-\lowscale})^{-1/2}\to0$.
Moreover, by Lemma~\ref{lem:etavareta},
$\sum_{(j,k)\in\indexset_n}(j-\jmean{\indexset_n}) j\sim
(n2^{-\lowscale}) (2-2^{-(\upscale-\lowscale)})\kappa_{\upscale
-\lowscale}$ as $n\to\infty$.
Hence,
\[
\sum_{(j,k)\in\indexset_n}(j-\jmean{\indexset_n}) j  \frac{\dwt
_{j,k}^2}{2^{2\td j}} = (n2^{-\lowscale}) \bigl\{\bigl(2-2^{-(\upscale
-\lowscale)}\bigr) \kappa_{\upscale-\lowscale}+o_\prob(1)\bigr\} ,
\]
and~(\ref{eq:htdVStd}) and the previous display yield
%
\begin{eqnarray}
\label{eq:AsympEquivStd}
&&\sum_{(j,k)\in\indexset_n}(j-\jmean{\indexset_n}) j  \frac{\dwt
_{j,k}^2}{2^{2\htd_n j}}\nonumber\\[-8pt]\\[-8pt]
&&\qquad = (n2^{-\lowscale})
\bigl\{\K \bigl(2-2^{-(\upscale-\lowscale)}\bigr) \kappa_{\upscale
-\lowscale
}+o_\prob(1)\bigr\} .\nonumber
\end{eqnarray}
Since $\kappa_\k>0$ for all $\k\geq1$ and $\kappa_\k\to2$ as $\k
\to
\infty$ (see Lemma~\ref{lem:etavareta}),
and since we assumed $\upscale-\lowscale\geq1$, the sequence
$(2-2^{-(\upscale-\lowscale)}) \kappa_{\upscale-\lowscale}$ is
bounded below by a positive
constant, so (\ref{eq:EstEqTaylor}),~(\ref{eq:ScltBound})
and~(\ref{eq:AsympEquivStd}) imply~(\ref{eq:var-bias}).
\end{pf*}
\begin{pf*}{Proof of Theorem~\protect\ref{thm:CLT}} 
Define $f^\ast(0)\eqdef\d\mes^\ast/\d\lambda_{|\lambda=0}$.
Since $\mes^\ast\in\calH(\beta,\gamma,\varepsilon)$ and $\mes
^\ast
(-\varepsilon,\varepsilon)>0$, we have $f^\ast(0)>0$.
Without loss of generality, we set $f^\ast(0)=~1$.
By Corollary~\ref{cor:C1gaussian}, conditions~(\ref{eq:BiasTermAssump})
and~(\ref{eq:rosenthalAssump}) hold with
$\K=\Kvar[\psi]{\td}$. Moreover, (\ref{eq:CondBiaisNegligeable})
implies that
$\lowscale^{-1}+\lowscale^2(n2^{-\lowscale})^{-1/4}\to0$,
so we may apply~(\ref{eq:AsympEquivStd}), which, with~(\ref{eq:EstEqTaylor}), gives
\begin{equation}\label{eq:rateDecompUniv}
(n2^{-\lowscale})^{1/2} (\hd_n-\td) =
\frac{(n2^{-\lowscale})^{-1/2} \Sclt_n(\td)}{2 \log(2) \K
(2-2^{-(\upscale-\lowscale)})  \kappa_{\upscale-\lowscale}}
\bigl(1+o_{\prob}(1)\bigr).
\end{equation}
Define $\Scltbar_n$ as $\Sclt_n$ in~(\ref{eq:ScltDef}), but with the
wavelet coefficients $\overline{\dwt}_{j,k}$ defined in
Corollary~\ref{cor:C1gaussian} replacing the wavelet coefficients
$\dwt
_{j,k}$. Let us write
\begin{equation}\label{eq:linearPartDecomp}
\hspace*{6pt}\Sclt_n(\td)= \bigl(\Sclt_n(\td)-\Scltbar_n(\td)\bigr) + \PE_f
[\Scltbar_n(\td)]
+ \bigl(\Scltbar_n(\td)-\PE_f[\Scltbar_n(\td)]
\bigr) .
\end{equation}
By Corollary~\ref{cor:C1gaussian}, using Minkowski's and Markov's
inequalities,~(\ref{eq:etaJ0}), $\nj{j}\leq n2^{-j}$ and
$\td+\alpha>(1+\beta)/2$, we obtain, as $n\to\infty$,
\[
\Sclt_n(\td)-\Scltbar_n(\td)=o_\prob( (n2^{-\lowscale})^{1/2}).
\]
Since $\sum_{(j,k)\in\indexset_n} (j-\jmean{\indexset_n}) =0$ and
$\PE
_f [ \dwt_{j,k}^2 ]=\sigma^{2}_{j}(\nu)$, we may write
\begin{eqnarray*}
\PE_f [ \Sclt_n(\td) ]&=& \sum_{(j,k)\in\indexset_n}
(j-\jmean
{\indexset_n})   \bigl(2^{-2\td j}\sigma^{2}_{j}(\nu)-\K\bigr)\\
&=& O\bigl(n2^{-(1+\beta) \lowscale}\bigr)
= o( (n2^{-\lowscale})^{1/2} ),
\end{eqnarray*}
where the $O$-term follows from~(\ref{eq:BiasTermAssump}), (\ref{eq:etaJ0}) and $\nj{j}\leq
n2^{-j}$ and the $o$-term follows from~(\ref{eq:CondBiaisNegligeable}).
Using~(\ref{eq:rateDecompUniv}),~(\ref{eq:linearPartDecomp}) and the
two last displays,
we finally get that
\[
(n2^{-\lowscale})^{1/2} (\hd_n-\td) = \frac{(n2^{-\lowscale
})^{-1/2}
(\Scltbar_n(\td)-\PE_f[\Scltbar_n(\td)]
)}{2 \log
(2) \K (2-2^{-(\upscale-\lowscale)})
 \kappa_{\upscale-\lowscale}}  \bigl(1+o_{\prob}(1)\bigr).
\]
Because $\overline{f}(\lambda)=|1-\rme^{-\rmi\lambda}|^{-2\td
}[f^\ast\1
_{[-\epsilon,\epsilon]}](\lambda)$ and
$f^\ast\1_{[-\epsilon,\epsilon]}\in\calH(\beta,\gamma',\pi)$
for some
$\gamma'>0$, we may apply Proposition~\ref{prop:CLTLinearPart} below to
determine the asymptotic behavior of $\Scltbar_n(\td)-\PE_f[\Scltbar
_n(\td)]$ as $n\to\infty$. Since
$\K=f^\ast(0)\Kvar[\psi]{\td}$ (Theorem~\ref{thm:AssumpCons}), this
yields the result and
completes the proof.
\end{pf*}

The following proposition provides a CLT when the condition on $\mes
^\ast$ is global, namely
$\mes^\ast\in\calH(\beta,\gamma,\pi)$. It covers the
cases (\hyperref[it:caseJ1-J0fixed]{S-1}), where $\upscale-\lowscale\to\k<\infty$
and~\mbox{\textup{(\hyperref[it:caseJ1isJ]{S-2})}}, where $\upscale-\lowscale\to\infty$.
\begin{proposition}\label{prop:CLTLinearPart}
Let $X$ be a Gaussian process having generalized spectral measure~\textup{(\ref{eq:fmodele})}
with $\td\in\Rset$ and $\mes^\ast\in\calH(\beta,\gamma,\pi)$, with
$f^\ast(0)\eqdef\d\mes^\ast/\d\lambda_{|\lambda=0}>0$, where
$\gamma>0$ and $\beta\in(0,2]$.
Let $\lowscale$ and $\upscale$ be two sequences satisfying~\textup{(\ref{eq:J0J1n})}
and suppose that $\lowscale^{-1}+(n2^{-\lowscale
})^{-1}\to
0$ and
$\upscale-\lowscale\to\k\in\{1,2,\dots,\infty\}$ as $n\to\infty$.
Then, as $n\to\infty$,
%
\begin{equation}\label{eq:CLTSbar}
\frac{(n2^{-\lowscale})^{-1/2} \{\Sclt_{\indexset
_n(\lowscale
,\upscale)}(\td)-\PE_f[\Sclt_{\indexset_n(\lowscale,\upscale
)}(\td
)]\}}
{2 \log(2) f^\ast(0) \Kvar[\psi]{\td} (2-2^{-(\upscale
-\lowscale
)}) \kappa_{\upscale-\lowscale}}
\cl\calN(0,\AsympVarWWE[\psi]{\td,\k}) ,
\end{equation}
where $\kappa_k$ is defined in~\textup{(\ref{eq:eta-k})} and $\AsympVarWWE
[\psi]{\td,\k}$ in~\textup{(\ref{eq:varLimite-k})} for $\k<\infty$
and $\AsympVarWWE[\psi]{\td,\infty}$ in~\textup{(\ref{eq:varLimite})}.
\end{proposition}
\begin{pf*}{Proof}
We take $f^\ast(0)=1$, without loss of generality.
As $n\to\infty$, since $\upscale-\lowscale\to\k$, we have $\kappa
_{\upscale
-\lowscale}\to\kappa_\k$, by setting, in the special case where
$\k=\infty$, $\kappa_\infty=2$; see Lemma~\ref{lem:etavareta}. This
gives the deterministic limit of the denominator
in~(\ref{eq:CLTSbar}). The limit distribution of the numerator is
obtained by applying Lemma~\ref{lem:CLTgaussianQuadForm}
below.
Let $A_n$ and $\Gamma_n$ be the square matrices indexed by the pairs
$(j,k),   (j,k) \in\indexset_n\times\indexset_n$ (in
lexicographic order) and defined as follows:\looseness=1
\begin{enumerate}[(1)]
\item[(1)]$A_n$ is the diagonal matrix such that
$[A_n]_{(j,k),(j,k)}=(n2^{-\lowscale})^{-1/2}\mathrm{sign}(j-\jmean
{\indexset_n})$ for all $(j,k)\in\indexset_n$;
\item[(2)]$\Gamma_n$ is the covariance matrix of the vector
$\displaystyle[|j-\jmean{\indexset_n}|^{1/2} 2^{-\td j} \dwt
_{j,k}]_{(j,k)\in\indexset_n}$.
\end{enumerate}
Let $\rho(A)$ denote the spectral radius of the square matrix $A$,
that is, the maximum of the absolute value of its eigenvalues.
Of course, $\rho[A_n]=(n2^{-\lowscale})^{-1/2}$. Moreover,
$\rho[\Gamma_n]\leq\sum_{j=\lowscale}^{\upscale}\rho[\Gamma_{n,j}]$,
where $\Gamma_{n,j}$ is the covariance matrix of the vector
$[|j-\jmean{\indexset_n}|^{1/2} 2^{-\td j} \dwt_{j,k}
]_{k=0,\dots,\nj{j}-1}$.
Since $\{\dwt_{j,k}\}_{k\in\Zset}$ is a stationary time series, by
Lemma~\ref{lem:CovRadBound},
\[
\rho[\Gamma_{n,j}] \leq|j-\jmean{\indexset_n}|   2^{-2\td j}
2\pi
\sup_{\lambda\in(-\pi,\pi)} \bdens{j,0}{\lambda}{\mes}  .
\]
From~(\ref{eq:DjApprox}), since $\bdensasymp{0}{\cdot}{\td}$ is bounded
on $(-\pi,\pi)$, we get, for a constant
$C$ not depending on $n$,
$\rho[\Gamma_n]\leq C  \sum_{j=\lowscale}^{\upscale} |j-\jmean
{\indexset_n}| .$
By~(\ref{eq:etaJ0}), the latter sum is $O((\upscale-\lowscale)^2)$.
Hence, as $n\to\infty$, since $\upscale-\lowscale\leq\nmaxscale
-\lowscale=O(\log(n2^{-\lowscale}))$, we have
$\rho[A_n]\rho[\Gamma_n]=O((n2^{-\lowscale})^{-1/2}(\upscale
-\lowscale)^2)\to0$, so the conditions of
Lemma~\ref{lem:CLTgaussianQuadForm} are met, provided that
$(n2^{-\lowscale})^{-1}\PVar(\Sclt_n(\td))$ has a finite
limit.

To conclude the proof, we need to compute this limit.
In \cite{moulinesrouefftaqqu2007a}, Proposition~2, it is shown that
for all $\dj=0,1,\dots,$ as $j\to\infty$ and
$\nj{j}\to\infty$,
%
\begin{equation}\label{eq:c_ndefandconv}
c_n(j,\dj)\eqdef
2^{-4\td j}\nj{j-\dj}
\PCov(\hat{\sigma}_{j}^{2},
\hat{\sigma}_{j-\dj}{2}
)\to4\pi
 \intbdens[\psi]{\dj}{\td} ,
\end{equation}
where $\intbdens[\psi]{\dj}{d}$ is defined in~(\ref{eq:bDint}) and
$\hat{\sigma}_{j}^{2}\eqdef\frac1{\nj{j}}   \sum_{k=0}^{\nj{j}-1} \dwt_{j,k}^2$.
Since $\Sclt_n(\td)=\sum_{j=\lowscale}^{\upscale}(j-\jmean
{\indexset
_n})2^{-2j\td}\nj{j}\hat{\sigma}_{j}^{2}$, we obtain
\begin{eqnarray}\label{eq:AssympVar1}
&&(n2^{-\lowscale})^{-1}\PVar(\Sclt_n(\td))\nonumber\\
&&\qquad =\sum_{i=0}^{\upscale-\lowscale}(i+\lowscale-\jmean{\indexset
_n})^22^{-i} \frac{\nj{\lowscale+i}}{n2^{-(\lowscale+i)}}
c_n(\lowscale+i,0)\nonumber\\[-8pt]\\[-8pt]
&&\hspace*{33pt}{}+2 \sum_{i=1}^{\upscale-\lowscale}\sum_{\dj=1}^{i}(i+\lowscale
-\jmean
{\indexset_n})(i-\dj+\lowscale-\jmean{\indexset_n})\nonumber\\
&&\hspace*{86pt}{}\times2^{2\td\dj-i}
\frac{\nj{\lowscale+i}}{n2^{-(\lowscale+i)}}   c_n(\lowscale
+i,\dj).\nonumber
\end{eqnarray}
By the Cauchy--Schwarz inequality,~(\ref{eq:rosenthalAssumpBis}),
(\ref{eq:BiasTermAssump}) and $\nj{j-\dj}\asymp\nj{j}2^{-\dj}$ imply
that $|c_n(j,\dj)|\leq C2^{-2\td\dj+\dj/2}$, where $C$ is a positive constant.
Using this bound, (\ref{eq:etaJ0}) and $\nj{j}\leq n2^{-j}$ for
bounding the terms of the two series in the right-hand side
of~(\ref{eq:AssympVar1}) yields the following convergent series: $\sum
_{i=0}^{\infty}(i+1)^22^{-i}$ and
$\sum_{i=1}^{\infty}\sum_{\dj=1}^{i}(i+1)(i-\dj+1) 2^{-i+\dj/2}$. Using
the assumptions on $\upscale$ and $\lowscale$,\break we have
$\nj{\lowscale+i}\sim n2^{-(\lowscale+i)}$ for any
$i\geq0$ and by Lemma~\ref{lem:etavareta}, $\jmean{\indexset
_n}-\lowscale\to\eta_\k$ as $n\to\infty$.
Hence, by dominated convergence, (\ref{eq:AssympVar1})
and~(\ref{eq:c_ndefandconv}) finally give that, as $n\to\infty$,
$(n2^{-\lowscale})^{-1}\PVar(\Sclt_n(\td))$ converges to
%
\begin{equation}\label{eq:AssympVarFin}
\hspace*{12pt}4\pi\Biggl[
\intbdens[\psi]{0}{\td} \kappa_l(2-2^{-\k})
+ 2 \sum_{1\leq\dj\leq i \leq\k}(i-\eta_\k)(i-\eta_\k-\dj
)2^{2\td\dj
-i}  \intbdens[\psi]{\dj}{\td}\Biggr],
\end{equation}
where in the case $\k=\infty$, we have set $2^{-\infty}=0$, $\eta
_\infty
=1$ and $\kappa_\infty=2$.
Note that the above bound on $|c_n(j,\dj)|$ and~(\ref
{eq:c_ndefandconv}) imply that as $\dj\to\infty$,
%
\begin{equation}
\label{eq:iuBound}
\intbdens[\psi]{\dj}{\td}=O(2^{-2\td\dj+\dj/2}),
\end{equation}
which confirms that the series in~(\ref{eq:AssympVarFin}) is
convergent for
$\k=\infty$. Finally, dividing this variance by the squared limit of the
denominator in~(\ref{eq:CLTSbar}), we get the limit variance in
(\ref{eq:CLTSbar}), namely~(\ref{eq:varLimite-k})
and~(\ref{eq:varLimite}).
\end{pf*}

The following lemmas were used in the proof of Proposition~\ref{prop:CLTLinearPart}.
%
\begin{lemma}\label{lem:CovRadBound}
Let $\{\xi_\k, \k\in\Zset\}$ be a stationary process with spectral
density $g$ and let $\Gamma_n$ be the covariance matrix of $[\xi_1,
\dots, \xi_n]$. Then, $ \rho(\Gamma_n)\leq2\pi \|g\|_\infty$.
\end{lemma}
%
%
\begin{lemma}\label{lem:CLTgaussianQuadForm}
Let $\{\xi_{n}, n\geq1\}$ be a sequence of Gaussian
vectors with zero mean and covariance $\Gamma_n$.
Let $(A_n)_{n\geq1}$ be a sequence of deterministic symmetric
matrices 
such that
$ \lim_{n\to\infty} \PVar(\xi_n^T A_n \xi_n) =
\sigma^2 \in
[0,\infty). $
Assume that $\lim_{n\to\infty} [\rho(A_n)\rho(\Gamma
_n)]= 0$.
Then,
$\xi_n^T A_n \xi_n - \PE[\xi_n^T A_n \xi_n] \cl
\calN(0,\sigma^2) $.
\end{lemma}
\begin{pf*}{Proof}
The result is obvious if $\sigma=0$, hence we may assume $\sigma>0$.
Let \mbox{$n\geq1$}, $k_n$ be the rank of $\Gamma_n$ and $Q_n$ denote an
$n\times k_n$ full-rank matrix such that $Q_nQ_n^T=\Gamma_n$. Let
$\zeta_n\sim\calN(0,I_{k_n})$, where $I_k$ is the identity matrix of
size $k\times k$. Then, for any $k_n \times k_n$ unitary matrix
$U_n$, $U_n \zeta_n\sim\calN(0,I_{k_n})$ and hence $Q_n U_n \zeta_n$
has the same distribution as $\xi_n$. Moreover, since $A_n$ is
symmetric, so is $Q_n^TA_nQ_n$. Choose $U_n$ to be a unitary matrix
such that $\Lambda_n\eqdef U_n^T(Q_n^TA_nQ_n)U_n$ is a diagonal
matrix. Thus, $ \zeta_n^T\Lambda_n\zeta_n= (Q_nU_n\zeta_n)^T A_n
(Q_nU_n\zeta_n)$ has the same distribution as $\xi_n^T A_n \xi_n$.
Since $\Lambda_n$ is diagonal, $\zeta_n^T\Lambda_n\zeta_n$ is a sum
of independent r.v.'s of the form $\sum_{k=1}^{k_n} \lambda_{k,n}
\zeta_{k,n}^2$, where $(\zeta_{1,n}, \dots, \zeta_{k_n,n})$ are
independent centered unit-variance Gaussian r.v.'s and
$\lambda_{k,n}$ are the diagonal entries of $\Lambda_n$. Note that
$\sum_{k=1}^{k_n} \lambda_{k,n}=\PE[\xi_n^T A_n \xi_n]$.
To check the asymptotic normality, we verify that the Lindeberg
conditions hold for the sum of centered independent r.v.'s: $\xi_n^T
A_n \xi_n - \PE[\xi_n^T A_n \xi_n] = \sum_{k=1}^{k_n}
\lambda_{k,n} (\zeta_{k,n}^2-1)$. Under the stated assumptions,
\[
\sum_{k=1}^{k_n}
\lambda_{k,n}^2 \PE(\zeta_{k,n}^2-1)^2 = \PVar(\xi_n^T A_n
\xi_n) \to\sigma^2 \qquad\mbox{as $n \to\infty$}
\]
and $\rho(\Lambda_n)=\rho(Q_n^TA_nQ_n)\leq
\rho(A_n) \sup_{\|x\|=1}\|Q_nx\|^2=\rho(A_n) \rho(\Gamma_n)\to0
.$ Since $\rho(\Lambda_n) = \max_{1\leq k\leq
k_n}|\lambda_{k,n}|$, for all $\epsilon>0$,
\begin{eqnarray*}
&&\sum_{k=1}^{k_n} \lambda_{k,n}^2 \PE\bigl[(\zeta_{k,n}^2-1)^2
\1\bigl(|\lambda_{k,n} (\zeta_{k,n}^2-1)|\geq\epsilon\bigr)\bigr]\\
&&\qquad\leq\Biggl(\sum_{k=1}^{k_n} \lambda_{k,n}^2\Biggr)
\PE\bigl[\bigl(\zeta_{1,n}^2-1\bigr)^2 \1\bigl(\rho(\Lambda_n)
 |\zeta_{1,n}^2-1|\geq\epsilon\bigr)\bigr] \to0\qquad\mbox{as $n
\to\infty$}.
\end{eqnarray*}
Hence, the Lindeberg conditions hold provided $\sigma>0$.
\end{pf*}
\begin{lemma}
\label{lem:etavareta}
Let $p,\k\geq0$, $\eta_\k$ and
$\kappa_\k$ be defined as in~\textup{(\ref{eq:eta-k})}, $\jmean{\indexset}$ as
in~\textup{(\ref{eq:etaDelta})} and
\[
\vareta(\indexset)\eqdef\cardinal{\indexset}^{-1}\sum_{(j,k)\in
\indexset}
(j-\jmean{\indexset})^2 =\cardinal{\indexset}^{-1}\sum_{(j,k)\in
\indexset}
j(j-\jmean{\indexset}) .
\]
We have
\begin{eqnarray}
&&\eta_\k= \frac{1 -2^{-\k}(1+\k/2)}{1-2^{-(\k+1)}} \in(0,1),\qquad
\lim_{\l\to\infty} \eta_\k= 1,\qquad\lim_{\l\to\infty} \kappa
_\k= 2,\label{eq:varetaLim}\\
&&\mbox{for all $\dj\geq0$}\qquad
\lim_{\k\to\infty}\frac1{\kappa_\k}\sum_{i=0}^{\k-\dj}\frac
{2^{-i}}{2-2^{-\k}}(i-\eta_\k)(i+\dj-\eta_\k)=1\label{eq:Varkinfty}
\end{eqnarray}
and for all $n\geq1$ and $0\leq j_0 \leq j_1\leq\nmaxscale$,
\begin{equation}
\label{eq:etavaretafirst}
\Biggl| \sum_{j=j_0}^{j_1} (j-j_0)^p \nj{j} -
n2^{-j_0} \sum_{i=0}^{j_1-j_0} i^p2^{-i}\Biggr| \leq2(\L-1)
(j_1-j_0)^{p+1}.
\end{equation}
Moreover, if $0\leq\nlowscale\leq\nmaxscale$ with
$n2^{-\nlowscale}\to\infty$ as $n\to\infty$, then
\begin{eqnarray*}
 \sup_{j_1=\nlowscale,\dots,\nmaxscale} \bigl|\cardinal
{\indexset
_n(\nlowscale,j_1)} - n2^{-\nlowscale}\bigl(2-2^{-(j_1-\nlowscale
)}\bigr)\bigr|&=& O(\log(n2^{-\nlowscale})) , \\
 \sup_{j_1=\nlowscale,\dots,\nmaxscale} |\jmean{\indexset
_n(\nlowscale,j_1)} - \nlowscale- \eta_{j_1-\nlowscale}| &=& O
(\log^2(n2^{-\nlowscale}) (n2^{-\nlowscale})^{-1}),\\
 \sup_{j_1=\nlowscale,\dots,\nmaxscale} |\vareta[\indexset
_n(\nlowscale,j_1)]-\kappa_{j_1-\nlowscale}| &=&
O(\log^3(n2^{-\nlowscale}) (n2^{-\nlowscale})^{-1}).
\end{eqnarray*}
\end{lemma}
%
%

\section*{Acknowledgments}

We would like to thank the referees
and the Associate Editor for their helpful comments.
Murad~S.~Taqqu
would like to
 thank l'\'Ecole Normale Sup\'erieure des T\'elecom\-munications in Paris
for
their hospitality.

\printaddresses

\end{document}